# Higher order temporal finite element methods through mixed formalisms

Jinkyu Kim

*Research Professor, School of Civil, Environmental and Architectural Engineering, Korea University, Anam-dong5-ga1, Seongbuk-goo 136-713, Korea.*

Phone: 82 2 3290 3833

Fax: 82 2 921 2439

Email: jk295@korea.ac.kr

## Abstract

The extended framework of Hamilton's principle and the mixed convolved action principle provide new rigorous weak variational formalism for a broad range of initial boundary value problems in mathematical physics and mechanics. Both approaches utilize the mixed formulation and lead to the development of various space-time finite element methods. In this paper, their potential when adopting temporally higher order approximations is investigated. The classical single-degree-of-freedom dynamical systems are primarily considered to validate and to investigate the performance of the numerical algorithms developed from both formulations. For the undamped system, all the algorithms are found to be symplectic and unconditionally stable with respect to the time step. On the other hand, for the damped system, the approach is shown to be robust and to be accurate with good convergence characteristics.

## 1. Introduction

Despite of its origin in particle dynamics, Hamilton's principle [1, 2] has been with us for a long time throughout broad range of mathematical physics [3-7]. However, it suffers from two main difficulties such as (i) use of end-point constraints and (ii) adoption of Rayleigh's dissipation for non-conservative systems. The first difficulty relates to the proper use of initial conditions resulting from the restrictions on the function variations. In Hamilton's principle, the variations vanish at the end points of the time interval, which, in turn, implies that the functions are known at these two instants. For a typical dynamic problem, one does not know how the considered system evolves at the end of the time interval. Usually, this is the main objective of the analysis, which means that there may be a serious philosophical or mathematical inconsistency in Hamilton's principle. Second difficulty relates to the inability to incorporate irreversible phenomena. Hamilton's principle itself only applies to conservative systems. With Rayleigh's dissipation [8], irreversible processes can be brought into the framework of Hamilton's principle. However, this approach is not satisfactory in a strict mathematical sense, since the variation of Rayleigh's dissipation enters in an ad-hoc manner.

Historically, to resolve such difficulties in Hamilton's principle, Tonti [9, 10] suggested that convolution should replace the inner product for variational methods in initial value problems. Somewhat earlier, Gurtin [11-13] introduced the convolution functional, and could reduce the initial value problem to an equivalent boundary value problem. However, the functional by Gurtin is complicated and it never can recover the original strong form. Following the ideas of Tonti and Gurtin, Oden and Reddy [14] extended the formulation to a large class of initial boundary problems in mechanics, especially for Hellinger-Reissner type mixed principles. More recently, Riewer [15, 16] adopted the use of fractional calculus to accommodate dissipative dynamical systems. This is an attractive idea, and many other researches including Agrawal [17-19], Atanackovic et al. [20], Baleanu and Muslih [21], Dreisigmmeyer and Young [22, 23], El-Nabusi and Torres [24], and Abreu and Godinho [25] have proposed similar approaches. However, surprisingly, none of these papers include an analytical description validating their approach for the most fundamental case, a classical Kelvin-Voigt single-degree-of-freedom (SDOF) damped oscillator.

Recently, two new variational frameworks for elastodynamics such as extended framework of Hamilton's principle (EHP, [26]) and mixed convolved action principle (MCAP, [27]) were established by using mixed variables. While EHP adopts a mixed Lagrangian formalism given in [28-31], it provides a new and simple framework that correctly accounts for initial conditions within Hamilton's principle. EHP resides in an incomplete variational framework since it requires Rayleigh's function for dissipative systems and cannot define the functional action, explicitly. On the other hand, MCAP clearly resolves long-standing problems in Hamilton's principle. With MCAP, a single scalar functional action provides the governing differential equations, along with all the pertinent boundary and initial conditions for conservative and non-conservative linear systems. Thus, in theoretical aspects, MCAP is certainly preferred rather than EHP, however, there still remains a challenge for MCAP to have the generalized framework of other than linear problems. While EHP can be numerically implemented for viscoplasticity continuum dynamics, MCAP is currently suffered to have the explicit functional action for such problems. Since both methods provide sound basis to develop various space-time finite element methods for linear initial boundary value problems, here the focus is initially on investigating their potential when employing higher-order temporal approximations.

The remainder of the paper is organized as follows. Next, in Section 2, we provide some relevant background on EHP and MCAP, especially for the SDOF Kelvin-Voigt system. In Section 3, discretization scheme and numerical algorithms are provided when temporally higher-order approximations are adopted in both approaches. Basic numerical properties of the developed methods are closely examined in Section 4. Then, some numerical examples are presented to investigate and to validate all of these developed algorithms for practical problems of the forced vibration in Section 5. Finally, some conclusions are provided in Section 6.

## 2. New variational formalisms

In this Section, new variational frameworks for the SDOF Kelvin-Voigt system displayed in Fig. 1 were reviewed for the development of higher order temporal finite element methods from both approaches.

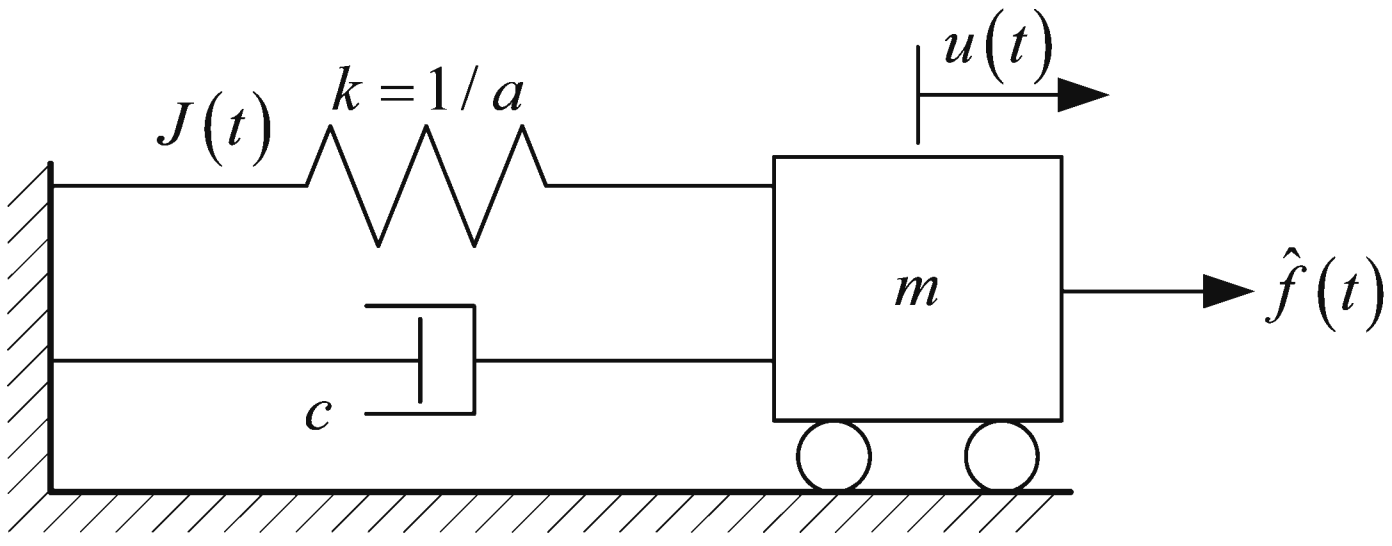

**Fig. 1.** SDOF Kelvin-Voigt damped oscillator

With mass $m$, damping coefficient $c$, the known applied force $\hat{f}(t)$ with time $t$, and stiffness $k=1/a$ with $a$ representing the flexibility, EHP and MCAP could formulate the variational framework for this model in terms of the displacement of the mass $u(t)$ and the impulse of the internal force $J(t)$ in the spring.

### *2.1. Weak form for the Kelvin-Voigt model in EHP*

Following the ideas in [31], the EHP associated with this problem defines Lagrangian $L$ and Rayleigh's dissipation $\varphi$ as

$$L(u,\dot{u},\dot{J};t)=\frac{1}{2}m\,\dot{u}^2+\frac{1}{2}a\,\dot{J}^2-\dot{J}\,u+\hat{f}\,u \tag{1}$$

and

$$\varphi(\dot{u};t)=\tfrac{1}{2}c\left[\dot{u}(t)\right]^2 \tag{2}$$

where a superposed dot represents a derivative with respect to time.

Then, the functional action $A$ for the fixed time interval from $t_0$ to $t$ is given by

$$A(u,\dot{u};t)=\int_{t_0}^{t} L(u,\dot{u};\tau)\,d\tau \tag{3}$$

and, in EHP, the first variation of $A$ is newly defined as

$$\delta A_{NEW} = -\delta \int_{t_0}^{t} L(u,\dot{u};\tau)d\tau + \int_{t_0}^{t} \frac{\partial \varphi(\dot{u};\tau)}{\partial \dot{u}} \delta u \, d\tau + \underline{\left[ m\hat{\dot{u}}\, \delta\hat{u} \right]_{t_0}^{t}} = 0 \tag{4}$$

by adding the counterparts (the underlined terms in Eq. (4)) to the terms without the end-point constraints in Hamilton's principle.

Such added terms have effect on confining a dynamical system to evolve uniquely from start to end with the unspecified values at the ends of the time interval such as $\hat{\dot{u}}(t_0)$, $\hat{u}(t_0)$, $\hat{\dot{u}}(t)$ and $\hat{u}(t)$. Then, interpreting the unspecified initial terms as sequentially assigning the known initial values completes this formulation. Thus, in EHP, the given initial velocity $\dot{u}_0$ is assigned first, and the given initial displacement $u_0$ is assigned next by

$$\hat{\dot{u}}(t_0) = \dot{u}_0 \tag{5}$$

and

$$\delta\hat{u}(t_0) = 0 \quad \left(\because u(t_0) = u_0\right) \tag{6}$$

The subsequent zero-valued term (6) needs not appear explicitly in the new action variation, so that the new definition (4) with the sequential assigning process such as (5) and (6) can properly account for the initial value problems. It should be noted that in EHP, the dependent initial condition $J_0$ can be identified by

$$m\dot{u}_0 + cu_0 + J_0 - \hat{j}_0 = 0 \tag{7}$$

where $\hat{j}_0$ is the initial internal impulse of the known applied force $\hat{f}$ given by

$$\hat{j}_0 = \int_{-\infty}^{t_0} \hat{f}(\tau)\, d\tau \tag{8}$$

In Eq. (8), the time interval $[-\infty, t_0]$ is used to represent that this is the time interval before the initial time we are considering.

To check this, let us substitute Eqs. (1)-(2) into Eq. (4). Then, we have

$$\delta A_{NEW} = -\int_{t_0}^{t} \left[ m\,\dot{u}\,\delta\dot{u} - c\,\dot{u}\,\delta u - \dot{J}\,\delta u + \hat{f}\,\delta u + a\,\dot{J}\,\delta\dot{J} - u\,\delta\dot{J} \right] d\tau + \left[ m\hat{\dot{u}}\,\delta\hat{u} \right]_{t_0}^{t} = 0 \tag{9}$$

Doing integration by parts on $m\dot{u}\,\delta\dot{u}$, $a\,\dot{J}\,\delta\dot{J}$, and $-u\,\delta\dot{J}$ in Eq. (9) yields

$$\begin{aligned}\delta A_{NEW} = & \underline{\left[ m\hat{\dot{u}}(t)\,\delta\hat{u}(t) - m\hat{\dot{u}}(t_0)\,\delta\hat{u}(t_0) \right]} - \left[ m\,\dot{u}\,\delta u \right]_{t_0}^{t} + \left[ \left( u - a\,\dot{J} \right) \delta J \right]_{t_0}^{t} \\ & + \int_{t_0}^{t} \left( m\,\ddot{u} + c\dot{u} + \dot{J} - \hat{f} \right) \delta u \, d\tau + \int_{t_0}^{t} \left( a\,\ddot{J} - \dot{u} \right) \delta J \, d\tau = 0 \end{aligned} \tag{10}$$

For arbitrary variations of $\delta u$ and $\delta J$ for the time interval $[t_0, t]$, the governing differential equations are given by

$$m\,\ddot{u}+c\,\dot{u}+\dot{J}-\hat{f}=0; \quad a\,\ddot{J}-\dot{u}=0 \tag{11}$$

along with constitutive relation as

$$u-a\,\dot{J}=0 \tag{12}$$

With the underlined terms in Eq. (10), the trajectory of the damped oscillator is firstly uniquely confined by

$$\hat{u}(t)=\dot{u}(t); \quad \delta\hat{u}(t)=\delta u(t) \tag{13}$$

while the given initial conditions are identified sequentially by Eq. (5), Eq. (6) and Eq. (7). Thus, with EHP, Hamilton's principle can account for compatible initial conditions to the strong form. It is not a complete variational method, since it still requires the Rayleigh's dissipation for a non-conservative process and the first variation of the functional action cannot yield the proper weak form explicitly. However, the framework is quite simple and it can be readily applied to problems other than linear elasticity with the use of Rayleigh's dissipation.

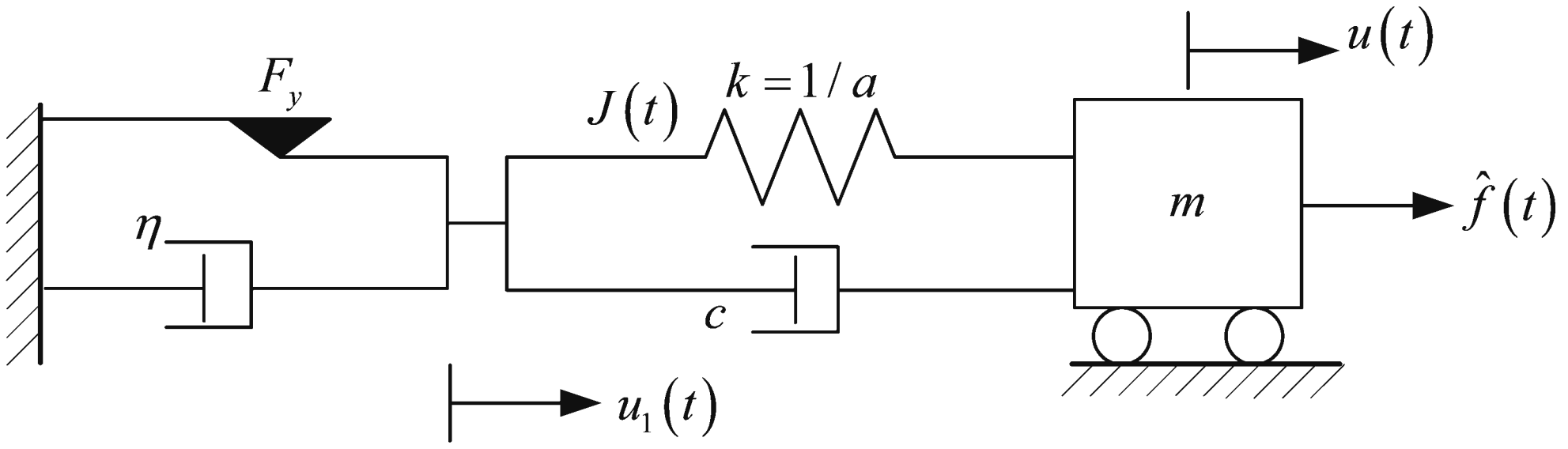

**Fig. 2.** SDOF elasto-viscoplastic model

For a representative example, let us consider SDOF elasto-viscoplastic model in Fig. 2. Rayleigh's dissipation to define rate-deformation for the slider-dashpot $\dot{u}_1$ can be given by

$$\varphi\left(\dot{J};\, t\right)=\frac{1}{2\eta}\left\langle\left|\dot{J}\right|-F_y\right\rangle^2 \tag{14}$$

in terms of Macaulay bracket $\langle\bullet\rangle$ and absolute value of $\dot{J}$ whereby $\eta$ and $F_y$ represent viscosity and yield force, respectively. Thus, in EHP, the action variation for this model is defined by adding up Eq. (4) and $\delta A_\varphi$

$$\begin{aligned}\delta A_\varphi &= \int_{t_0}^{t}\frac{\partial\varphi\left(\dot{J};\tau\right)}{\partial\dot{J}}\delta J\,d\tau \\ &= \int_{t_0}^{t}\underline{\frac{1}{\eta}\left\langle\left|\dot{J}\right|-F_y\right\rangle\operatorname{sgn}\left(\dot{J}\right)}\delta J\,d\tau\end{aligned} \tag{15}$$

and

$$-\left[\hat{u}_1\,\delta\hat{J}\right]_{t_0}^{t} \tag{16}$$

where the underlined term represents the rate-deformation for the slider-dashpot, $\dot{u}_1$.

Note that the adding terms (16) are the counterparts to the terms without the end-point constraints in Hamilton's principle that obtained from the compatibility condition

$$a\dot{J} + u_1 = u \tag{17}$$

With Eq. (4) and Eqs. (14)-(16), the governing differential equations for Fig. 2

$$m\ddot{u} + c\dot{u} + \dot{J} = \hat{f}; \quad a\ddot{J} - \dot{u} + \dot{u}_1 = 0 \tag{18}$$

are properly recovered in EHP along with proper initial conditions such as Eqs. (5)-(7) and $\hat{u}_1$ at $t = 0$.

### *2.2. Weak form for the Kelvin-Voigt model in MCAP*

As well described in [27], MCAP defines the convolved action for the SDOF Kelvin-Voigt damped oscillator as

$$A(u, \breve{u}, \dot{u}, J, \breve{J}, \dot{J}; t) = \frac{1}{2}m(\dot{u}*\dot{u}) - \frac{1}{2}a(\dot{J}*\dot{J}) + (\breve{J}*\breve{u}) + \frac{1}{2}c(\breve{u}*\breve{u}) - (u*\hat{f}) - u(t)\,\hat{j}(0) \tag{19}$$

where a superimposed arc represents a temporal left Riemann-Liouville semi-derivative. Referred to [32, 33], this is defined by

$$\breve{u}(t) = \left(D_{0^+}^{1/2}u\right)(t) \equiv \frac{1}{\Gamma(1/2)}\frac{d}{dt}\int_0^t \frac{u(\tau)}{(t-\tau)^{1/2}}\,d\tau \tag{20}$$

where $\Gamma(\bullet)$ denotes the Gamma function.

In Eq. (19), the symbol $*$ represents the convolution of two functions over time, such that

$$(\varphi*\phi)(t) = \int_0^t \varphi(\tau)\phi(t-\tau)\,d\tau \tag{21}$$

Meanwhile, the last term $\hat{j}(0)$ in Eq. (19) represents the initial impulse corresponding to $\hat{f}(t)$ that given by

$$\hat{j}(t) = \int_{-\infty}^{t} \hat{f}(\tau)\,d\tau \tag{22}$$

In MCAP, the stationarity of the action (19) yields the following weak form in time

$$\delta A = m(\delta\dot{u}*\dot{u}) - a(\delta\dot{J}*\dot{J}) + (\delta\breve{J}*\breve{u}) + (\delta\breve{u}*\breve{J}) + c(\delta\breve{u}*\breve{u}) - (\delta u*\hat{f}) - \delta u(t)\,\hat{j}(0) = 0 \tag{23}$$

After performing classical and fractional integration by parts on the appropriate terms in Eq. (23) as follows [27], we have

$$\begin{aligned}\delta A = &\left(\delta u*\left\{m\ddot{u} + c\dot{u} + \dot{J} - \hat{f}\right\}\right) + \left(\delta J*\left\{-a\ddot{J} + \dot{u}\right\}\right) + \delta u(t)\left\{m\dot{u}(0) + cu(0) + J(0) - \hat{j}(0)\right\} \\ &+ \delta u(0)\left\{m\dot{u}(t)\right\} + \delta J(t)\left\{-a\dot{J}(0) + u(0)\right\} - \delta J(0)\left\{-a\dot{J}(t)\right\} = 0\end{aligned} \tag{24}$$

For the sake of completeness, the fractional integration by parts formula is given

$$\int_0^t \left(D_{0^+}^{1/2}\varphi\right)(t)\ \left(D_{0^+}^{1/2}\phi\right)(t-\tau)\,d\tau = \int_0^t \frac{d\varphi}{d\tau}(\tau)\ \phi(t-\tau)\,d\tau + \varphi(0)\phi(t) \tag{25}$$

For arbitrary variations of $u$ and $J$, Eq. (24) emanates the governing differential equations in mixed forms as

$$m\,\ddot{u} + c\dot{u} + \dot{J} - \hat{f} = 0; \quad a\,\ddot{J} - \dot{u} = 0 \tag{26}$$

along with the proper initial conditions

$$m\dot{u}(0) + cu(0) + J(0) - \hat{j}(0) = 0; \quad -a\,\dot{J}(0) + u(0) = 0 \tag{27}$$

Note that the initial variations such as $\delta u(0)$ and $\delta J(0)$ vanish due to Eq. (27). In other words, in MCAP, we can identify the dependent initial conditions such as $J(0)$ and $\hat{j}(0)$ from the usual given initial conditions $u(0)$ and $\dot{u}(0)$ as well as the known initial impulse $\hat{j}(0)$.

As shown in Eq. (24) and Eq. (26)-(27), every governing equations and initial conditions are satisfied weakly in MCAP, where it incorporates both conservative and non-conservative components within the unified functional action (19). Thus, it resolves the long-standing problem in Hamilton's principle. However, MCAP still requires a generalized framework to embrace various irreversible phenomena. In particular, currently, it does not have the functional action for the problem shown in Fig. 2. Also, it should be noted that any pair of complementary order of fractional derivatives in Eq. (19) yields Eqs. (26)-(27) due to properties of complementary order of fractional derivatives in the fractional integration by parts

$$\int_0^t \left(D_{0^+}^{1-\alpha}\varphi\right)(t)\ \left(D_{0^+}^{\alpha}\phi\right)(t-\tau)\,d\tau = \int_0^t \frac{d\varphi}{d\tau}(\tau)\ \phi(t-\tau)\,d\tau + \varphi(0)\phi(t) \tag{28}$$

for $0 < \alpha < 1$.

## 3. Numerical implementation

The weak form (9) in EHP and the weak form (23) in MCAP include, at most, first derivatives of the primary variables $u(t)$ and $J(t)$ as well as the variations $\delta u(t)$ and $\delta J(t)$. Consequently, we have $C^0$ temporal continuity requirement on primary variables and the variations, thus, there are many cases to develop higher order temporal finite element methods. As we shall see in this Section, three kinds of quadratic temporal finite element methods in each framework are developed, since they are practically sufficient and accurate in computational aspects as discussed next. The numerical methods developed here are summarized in Table 1.

**Table 1** Developed quadratic temporal finite element methods in each framework

| Algorithms | Description |
|---|---|
| Jquad | $J(t)$ and $\delta J(t)$: quadratically approximated.<br>$u(t)$ and $\delta u(t)$ : linearly approximated |
| Uquad | $u(t)$ and $\delta u(t)$: quadratically approximated.<br>$J(t)$ and $\delta J(t)$: linearly approximated |
| UJquad | $u(t)$ and $\delta u(t)$: quadratically approximated.<br>$J(t)$ and $\delta J(t)$: quadratically approximated |

### *3.1. Algorithms from EHP*

By introducing the fixed time step $h$ for each time duration, that is, $t_r = r\,h$, Eq. (9) can be written

$$\delta A = \sum_{r=1}^{N} \delta A_r$$
$$= \sum_{r=1}^{N} -\left( \int_{t_{r-1}}^{t_r} \left[ m\,\dot{u}\,\delta\dot{u} - c\,\dot{u}\,\delta u - \dot{J}\,\delta u + \hat{f}\,\delta u + a\,\dot{J}\,\delta\dot{J} - u\,\delta\dot{J} \right] d\tau - \left[ \hat{p}\,\delta\hat{u} \right]_{t_{r-1}}^{t_r} \right) = 0 \tag{29}$$

where $\delta A_r$ represents the action variation in the $r^{th}$ time duration $[t_{r-1}, t_r]$ and $t_N = t$.
For $t_{r-1} \le \tau \le t_r$, temporally linear shape functions such as $L_{r-1}$ at $t_{r-1}$ and $L_r$ at $t_r$ are given by

$$L_{r-1}(\tau) = \frac{1}{h}(t_r - \tau) \tag{30}$$

$$L_r(\tau) = \frac{1}{h}(\tau - t_{r-1}) \tag{31}$$

Also, by introducing the center point $t_c$ for the time interval $[t_{r-1},\, t_r]$ as

$$t_c = \frac{(t_r - t_{r-1})}{2} = \frac{h}{2} \tag{32}$$

temporally quadratic shape functions $Q_{r-1}$ at $t_{r-1}$, $Q_r$ at $t_r$, and $Q_c$ at $t_c$ can be written as

$$Q_{r-1}(\tau)=\frac{2}{h^2}(\tau-t_r)(\tau-t_c) \tag{33}$$

$$Q_r(\tau)=\frac{2}{h^2}(\tau-t_{r-1})(\tau-t_c) \tag{34}$$

$$Q_c(\tau)=-\frac{4}{h^2}(\tau-t_r)(\tau-t_{r-1}) \tag{35}$$

With linear temporal shape functions (30)-(31) and quadratic temporal shape functions (33)-(35), we can develop every algorithms of EHP presented in Table 1.

For a representative case, Jquad algorithm can be obtained from the main approximations as

$$J(\tau)=Q_{r-1}(\tau)J_{r-1}+Q_r(\tau)J_r+Q_c(\tau)J_c \tag{36}$$

$$\delta J(\tau)=Q_{r-1}(\tau)\delta J_{r-1}+Q_r(\tau)\delta J_r+Q_c(\tau)\delta J_c \tag{37}$$

$$u(\tau)=L_{r-1}(\tau)u_{r-1}+L_r(\tau)u_r \tag{38}$$

$$\delta u(\tau)=L_{r-1}(\tau)\delta u_{r-1}+L_r(\tau)\delta u_r \tag{39}$$

$$\hat{f}(\tau)=L_{r-1}(\tau)\hat{f}_{r-1}+L_r(\tau)\hat{f}_r \tag{40}$$

and the subsequent approximations as

$$\dot{J}(\tau)=\dot{Q}_{r-1}(\tau)J_{r-1}+\dot{Q}_r(\tau)J_r+\dot{Q}_c(\tau)J_c \tag{41}$$

$$\delta\dot{J}(\tau)=\dot{Q}_{r-1}(\tau)\delta J_{r-1}+\dot{Q}_r(\tau)\delta J_r+\dot{Q}_c(\tau)\delta J_c \tag{42}$$

$$\dot{u}(\tau)=\dot{L}_{r-1}(\tau)u_{r-1}+\dot{L}_r(\tau)u_r \tag{43}$$

$$\delta\dot{u}(\tau)=\dot{L}_{r-1}(\tau)\delta u_{r-1}+\dot{L}_r(\tau)\delta u_r \tag{44}$$

Substituting Eqs. (36)-(44) into Eq. (29), and integrating yields

$$\delta A_r = \left( \frac{m}{h}\{u_r - u_{r-1}\} + \frac{c}{2}\{u_r - u_{r-1}\} - \hat{p}_{r-1} - \frac{5}{6}J_{r-1} + \frac{1}{6}J_r + \frac{2}{3}J_c - \frac{h}{3}\hat{f}_{r-1} - \frac{h}{6}\hat{f}_r \right)\delta u_{r-1}$$

$$+ \left( -\frac{m}{h}\{u_r - u_{r-1}\} + \frac{c}{2}\{u_r - u_{r-1}\} + \hat{p}_r - \frac{1}{6}J_{r-1} + \frac{5}{6}J_r - \frac{2}{3}J_c - \frac{h}{6}\hat{f}_{r-1} - \frac{h}{3}\hat{f}_r \right)\delta u_r$$

$$+ \left( -a\left\{ \frac{7}{3h}J_{r-1} + \frac{1}{3h}J_r - \frac{8}{3h}J_c \right\} - \frac{5}{6}u_{r-1} - \frac{1}{6}u_r \right)\delta J_{r-1} \tag{45}$$

$$+ \left( -a\left\{ \frac{1}{3h}J_{r-1} + \frac{7}{3h}J_r - \frac{8}{3h}J_c \right\} + \frac{1}{6}u_{r-1} + \frac{5}{6}u_r \right)\delta J_r$$

$$+ \underline{\left( -a\left\{ -\frac{8}{3h}J_{r-1} - \frac{8}{3h}J_r + \frac{16}{3h}J_c \right\} + \frac{2}{3}u_{r-1} - \frac{2}{3}u_r \right)}\delta J_c = 0$$

By making the coefficient of $(\delta u_{r-1}, \delta u_r, \delta J_{r-1}, \delta J_r, \delta J_c)$ equal to zero in Eq. (45), we have four independent equations given by

$$\left( \frac{m}{h}\{u_r - u_{r-1}\} + \frac{c}{2}\{u_r - u_{r-1}\} - \hat{p}_{r-1} - \frac{5}{6}J_{r-1} + \frac{1}{6}J_r + \frac{2}{3}J_c - \frac{h}{3}\hat{f}_{r-1} - \frac{h}{6}\hat{f}_r \right) = 0 \tag{46}$$

$$\left( -\frac{m}{h}\{u_r - u_{r-1}\} + \frac{c}{2}\{u_r - u_{r-1}\} + \hat{p}_r - \frac{1}{6}J_{r-1} + \frac{5}{6}J_r - \frac{2}{3}J_c - \frac{h}{6}\hat{f}_{r-1} - \frac{h}{3}\hat{f}_r \right) = 0 \tag{47}$$

$$\left( -a\left\{ \frac{7}{3h}J_{r-1} + \frac{1}{3h}J_r - \frac{8}{3h}J_c \right\} - \frac{5}{6}u_{r-1} - \frac{1}{6}u_r \right) = 0 \tag{48}$$

$$\left( -a\left\{ \frac{1}{3h}J_{r-1} + \frac{7}{3h}J_r - \frac{8}{3h}J_c \right\} + \frac{1}{6}u_{r-1} + \frac{5}{6}u_r \right) = 0 \tag{49}$$

While deriving Eqs. (46)-(49), the equation from the underlined term is discarded because it is not independent, which can be obtained from adding Eq. (48) and Eq. (49).
From either Eq. (48) or Eq. (49), we can express $J_c$ in terms of $J_{r-1}$, $J_r$, $u_{r-1}$, and $u_r$, which finally yields the matrix equation of Jquad algorithm as

$$\begin{bmatrix} \frac{1}{12}\frac{(X+6cha)}{ha} & 0 & \frac{1}{2} \\ \frac{1}{12}\frac{(-X+6cha)}{ha} & 1 & \frac{1}{2} \\ 1 & 0 & -\frac{2a}{h} \end{bmatrix} \begin{Bmatrix} u_r \\ \hat{p}_r \\ J_r \end{Bmatrix} = \begin{bmatrix} \frac{1}{12}\frac{(X+6cha)}{ha} & 1 & \frac{1}{2} \\ \frac{1}{12}\frac{(-X+6cha)}{ha} & 0 & \frac{1}{2} \\ -1 & 0 & -\frac{2a}{h} \end{bmatrix} \begin{Bmatrix} u_{r-1} \\ \hat{p}_{r-1} \\ J_{r-1} \end{Bmatrix} + \begin{Bmatrix} \frac{h}{3}\hat{f}_{r-1} + \frac{h}{6}\hat{f}_r \\ \frac{h}{6}\hat{f}_{r-1} + \frac{h}{3}\hat{f}_r \\ 0 \end{Bmatrix} \tag{50}$$

where $X$ is given by

$$X = 12ma - h^2 \tag{51}$$

Similarly, we have the Uquad algorithm as

$$\begin{bmatrix} -\frac{1}{3}\frac{(3m+ch)}{h} & 0 & \frac{1}{6}\frac{(Y+6cha)}{h^2} \\ \frac{1}{3}\frac{(-9m+2ch)}{h} & -m & -\frac{1}{6}\frac{(-Y+6cha)}{h^2} \\ \frac{2}{3}\frac{(6m+ch)}{h} & 0 & -\frac{2}{3}\frac{X}{h^2} \end{bmatrix} \begin{Bmatrix} u_r \\ \hat{p}_r \\ J_r \end{Bmatrix} = \begin{bmatrix} \frac{1}{3}\frac{(9m+2ch)}{h} & -m & \frac{1}{6}\frac{(Y+6cha)}{h^2} \\ -\frac{1}{3}\frac{(-3m+ch)}{h} & 0 & -\frac{1}{6}\frac{(-Y+6cha)}{h^2} \\ \frac{2}{3}\frac{(-6m+ch)}{h} & 0 & -\frac{2}{3}\frac{X}{h^2} \end{bmatrix} \begin{Bmatrix} u_{r-1} \\ \hat{p}_{r-1} \\ J_{r-1} \end{Bmatrix} + \begin{Bmatrix} \frac{h}{6}\hat{f}_{r-1} \\ \frac{h}{6}\hat{f}_r \\ \frac{h}{3}\hat{f}_{r-1}+\frac{h}{6}\hat{f}_r \end{Bmatrix} \tag{52}$$

where $Y$ is given by

$$Y = 24ma + h^2 \tag{53}$$

Also, we have the UJquad algorithm as

$$\begin{bmatrix} -\frac{1}{12}\frac{(12ma+4cha+h^2)}{ha} & 0 & \frac{1}{6}\frac{(Y+6cha)}{h^2} \\ \frac{1}{12}\frac{(-36ma+8cha+h^2)}{ha} & 1 & \frac{1}{6}\frac{(Y-6cha)}{h^2} \\ \frac{2}{3}\frac{(6m+ch)}{h} & 0 & \frac{2}{3}\frac{(-X)}{h^2} \end{bmatrix} \begin{Bmatrix} u_r \\ \hat{p}_r \\ J_r \end{Bmatrix} = \begin{bmatrix} \frac{1}{12}\frac{(36ma+8cha-h^2)}{ha} & 1 & \frac{1}{6}\frac{(Y+6cha)}{h^2} \\ -\frac{1}{12}\frac{(-12ma+4cha-h^2)}{ha} & 0 & \frac{1}{6}\frac{(Y-6cha)}{h^2} \\ \frac{2}{3}\frac{(-6m+ch)}{h} & 0 & \frac{2}{3}\frac{(-X)}{h^2} \end{bmatrix} \begin{Bmatrix} u_{r-1} \\ \hat{p}_{r-1} \\ J_{r-1} \end{Bmatrix} + \begin{Bmatrix} \frac{h}{6}\hat{f}_{r-1} \\ \frac{h}{6}\hat{f}_r \\ \frac{h}{3}\hat{f}_{r-1}+\frac{h}{6}\hat{f}_r \end{Bmatrix} \tag{54}$$

with the adequate substitution of $u_c$ and $J_c$ in terms of $J_{r-1}$, $J_r$, $u_{r-1}$, and $u_r$.

### *3.2. Algorithms from MCAP*

Previously, MCAP was numerically implemented through linear temporal shape functions for classical SDOF oscillators and systems that utilize fractional-derivative constitutive models by [34]. Here, continuing through this line, but, the quadratic temporal finite element methods are developed.
As well described in [34], for any non-negative integer $m$ and $n$, we have the following relation

$$\left[\left(D_{0^+}^{1/2} t^m\right) * \left(D_{0^+}^{1/2} t^n\right)\right](t) = \frac{\Gamma(1+m)\Gamma(1+n)}{\Gamma(1+m+n)}\, t^{m+n} \tag{55}$$

for the convolution of the semi-derivatives of power functions.
To evaluate the convolution of semi-derivatives of polynomial shape functions, here, Eq. (55) is frequently used.

Since we cannot have summation form of the action variation in convolution integral (that is, $\delta A \neq \sum_{r=1}^{N} \delta A_r$), let us consider the action variation over one time-step $[0,\, h]$ as

$$\begin{aligned} &\delta A(u,\, \breve{u},\, \dot{u},\, J,\, \breve{J}, \dot{J};\, \underline{t \to h}) = \\ &m(\delta\dot{u} * \dot{u}) - a\left(\delta\dot{J} * \dot{J}\right) + \left(\delta\breve{J} * \breve{u}\right) + \left(\delta\breve{u} * \breve{J}\right) + c(\delta\breve{u} * \breve{u}) - \left(\delta u * \hat{f}\right) - \delta u(t)\, \hat{j}(0) = 0 \end{aligned} \tag{56}$$

where temporally linear and quadratic shape functions of $t(0 \le t \le h)$ are defined as

$$L_0(t) = 1 - \frac{t}{h} \tag{57}$$

$$L_1(t) = \frac{t}{h} \tag{58}$$

$$Q_0(t) = \frac{2}{h^2}\left(t^2 - \frac{3}{2}ht + \frac{h^2}{2}\right) \tag{59}$$

$$Q_1(t) = \frac{2}{h^2}\left(t^2 - \frac{h}{2}t\right) \tag{60}$$

$$Q_c(t) = -\frac{4}{h^2}\left(t^2 - ht\right) \tag{61}$$

Then, subsequent approximations are given by

$$\dot{L}_0(t) = -\frac{1}{h} \tag{62}$$

$$\dot{L}_1(t) = \frac{1}{h} \tag{63}$$

$$\dot{Q}_0(t) = \frac{2}{h^2}\left(2t - \frac{3}{2}h\right) \tag{64}$$

$$\dot{Q}_1(t) = \frac{2}{h^2}\left(2t - \frac{h}{2}\right) \tag{65}$$

$$\dot{Q}_c(t) = -\frac{4}{h^2}(2t - h) \tag{66}$$

Now, let us consider Jquad algorithm for a representative one.
With approximations (57)-(66), the convolution component $\left(\delta\breve{J} * \breve{u}\right)$ in Eq. (56) can be written as

$$\begin{aligned}
\left(\delta\breve{J} * \breve{u}\right)(t) &= \lfloor \delta J_0 \quad \delta J_1 \quad \delta J_c \rfloor \begin{Bmatrix} \breve{Q}_0 \\ \breve{Q}_1 \\ \breve{Q}_c \end{Bmatrix} * \lfloor \breve{L}_0 \quad \breve{L}_1 \rfloor \begin{Bmatrix} u_0 \\ u_1 \end{Bmatrix} \\
&= \lfloor \delta J_0 \quad \delta J_1 \quad \delta J_c \rfloor \begin{bmatrix} \breve{Q}_0 * \breve{L}_0 & \breve{Q}_0 * \breve{L}_1 \\ \breve{Q}_1 * \breve{L}_0 & \breve{Q}_1 * \breve{L}_1 \\ \breve{Q}_c * \breve{L}_0 & \breve{Q}_c * \breve{L}_1 \end{bmatrix} \begin{Bmatrix} u_0 \\ u_1 \end{Bmatrix}
\end{aligned} \tag{67}$$

in terms of row vector $\lfloor \bullet \rfloor$, matrix $[\bullet]$, and column vector $\{\bullet\}$.
Each component of matrix in Eq. (67) can be directly evaluated by using Eq. (55). For a representative one, $\left(\breve{Q}_0 * \breve{L}_0\right)(t)$ is computed as

$$\begin{aligned}
\left(\breve{Q}_0 * \breve{L}_0\right)(t) &= \left[\left(D_{0^+}^{1/2}\left\{\frac{2}{h^2}\left(t^2 - \frac{3}{2}ht + \frac{h^2}{2}\right)\right\}\right) * \left(D_{0^+}^{1/2}\left\{1 - \frac{t}{h}\right\}\right)\right](t) \\
&= \frac{2}{h^2}\left[t^2 - \frac{1}{h}\frac{t^3}{3} - \frac{3}{2}ht + \frac{3}{2}\frac{t^2}{2} + \frac{h^2}{2} - \frac{h}{2}t\right]
\end{aligned} \tag{68}$$

Then, by letting $t \to h$ in Eq. (68) due to the underlined term in Eq. (56), Eq. (68) yields

$$\begin{aligned}
\left(\breve{Q}_0 * \breve{L}_0\right)(t \to h) &= \frac{2}{h^2}\left[t^2 - \frac{1}{h}\frac{t^3}{3} - \frac{3}{2}ht + \frac{3}{2}\frac{t^2}{2} + \frac{h^2}{2} - \frac{h}{2}t\right] \\
&= \frac{2}{h^2}\left[h^2 - \frac{1}{h}\frac{h^3}{3} - \frac{3}{2}hh + \frac{3}{2}\frac{h^2}{2} + \frac{h^2}{2} - \frac{h}{2}h\right] \\
&= -\frac{1}{6}
\end{aligned} \tag{69}$$

Following the same procedures as in Eqs. (68)-(69), one finds

$$\delta \breve{J} * \breve{u} = \lfloor \delta J_0 \quad \delta J_1 \quad \delta J_c \rfloor \begin{bmatrix} -\frac{1}{6} & \frac{1}{6} \\ \frac{5}{6} & \frac{1}{6} \\ -\frac{2}{3} & \frac{2}{3} \end{bmatrix} \begin{Bmatrix} u_0 \\ u_1 \end{Bmatrix} \tag{70}$$

In a similar way,

$$\delta \breve{u} * \breve{J} = \lfloor \delta u_0 \quad \delta u_1 \rfloor \begin{bmatrix} -\frac{1}{6} & \frac{5}{6} & -\frac{2}{3} \\ \frac{1}{6} & \frac{1}{6} & \frac{2}{3} \end{bmatrix} \begin{Bmatrix} J_0 \\ J_1 \\ J_c \end{Bmatrix} \tag{71}$$

and for the viscous dissipation term

$$c\left(\delta \breve{u} * \breve{u}\right) = \lfloor \delta u_0 \quad \delta u_1 \rfloor \begin{bmatrix} -\frac{c}{2} & \frac{c}{2} \\ \frac{c}{2} & \frac{c}{2} \end{bmatrix} \begin{Bmatrix} u_0 \\ u_1 \end{Bmatrix} \tag{72}$$

With evaluation of typical integer order convolution components in Eq. (56), we have the following discretized weak form of Jquad:

$$\begin{aligned} &\lfloor \delta u_0 \quad \delta u_1 \rfloor \begin{bmatrix} \frac{m}{h} & -\frac{m}{h} \\ -\frac{m}{h} & \frac{m}{h} \end{bmatrix} \begin{Bmatrix} u_0 \\ u_1 \end{Bmatrix} - \lfloor \delta J_0 \quad \delta J_1 \quad \delta J_c \rfloor \begin{bmatrix} -\frac{a}{3h} & -\frac{7a}{3h} & \frac{8a}{3h} \\ -\frac{7a}{3h} & -\frac{a}{3h} & \frac{8a}{3h} \\ \frac{8a}{3h} & \frac{8a}{3h} & -\frac{16a}{3h} \end{bmatrix} \begin{Bmatrix} J_0 \\ J_1 \\ J_c \end{Bmatrix} \\ &+ \lfloor \delta J_0 \quad \delta J_1 \quad \delta J_c \rfloor \begin{bmatrix} -\frac{1}{6} & \frac{1}{6} \\ \frac{5}{6} & \frac{1}{6} \\ -\frac{2}{3} & \frac{2}{3} \end{bmatrix} \begin{Bmatrix} u_0 \\ u_1 \end{Bmatrix} + \lfloor \delta u_0 \quad \delta u_1 \rfloor \begin{bmatrix} -\frac{1}{6} & \frac{5}{6} & -\frac{2}{3} \\ \frac{1}{6} & \frac{1}{6} & \frac{2}{3} \end{bmatrix} \begin{Bmatrix} J_0 \\ J_1 \\ J_c \end{Bmatrix} \\ &+ \lfloor \delta u_0 \quad \delta u_1 \rfloor \begin{bmatrix} -\frac{c}{2} & \frac{c}{2} \\ \frac{c}{2} & \frac{c}{2} \end{bmatrix} \begin{Bmatrix} u_0 \\ u_1 \end{Bmatrix} - \lfloor \delta u_0 \quad \delta u_1 \rfloor \begin{bmatrix} \frac{h}{6} & \frac{h}{3} \\ \frac{h}{3} & \frac{h}{6} \end{bmatrix} \begin{Bmatrix} \hat{f}_0 \\ \hat{f}_1 \end{Bmatrix} - \delta u_1 \left[ \hat{j}_0 \right] = 0 \end{aligned} \tag{73}$$

With the known initial conditions $u_0$ and $J_0$, the variations $\delta u_0$ and $\delta J_0$ vanish. Thus, the weak form reduces to the following:

$$\delta u_1 \left\lfloor -\frac{m}{h} \quad \frac{m}{h} \right\rfloor \begin{Bmatrix} u_0 \\ u_1 \end{Bmatrix} + \lfloor \delta J_1 \quad \delta J_c \rfloor \begin{bmatrix} \frac{7a}{3h} & \frac{a}{3h} & -\frac{8a}{3h} \\ -\frac{8a}{3h} & -\frac{8a}{3h} & \frac{16a}{3h} \end{bmatrix} \begin{Bmatrix} J_0 \\ J_1 \\ J_c \end{Bmatrix}$$

$$+\lfloor \delta J_1 \quad \delta J_c \rfloor \begin{bmatrix} \frac{5}{6} & \frac{1}{6} \\ -\frac{2}{3} & \frac{2}{3} \end{bmatrix} \begin{Bmatrix} u_0 \\ u_1 \end{Bmatrix} + \delta u_1 \left\lfloor \frac{1}{6} \quad \frac{1}{6} \quad \frac{2}{3} \right\rfloor \begin{Bmatrix} J_0 \\ J_1 \\ J_c \end{Bmatrix} \tag{74}$$

$$+\delta u_1 \left\lfloor \frac{c}{2} \quad \frac{c}{2} \right\rfloor \begin{Bmatrix} u_0 \\ u_1 \end{Bmatrix} - \delta u_1 \left\lfloor \frac{h}{3} \quad \frac{h}{6} \right\rfloor \begin{Bmatrix} \hat{f}_0 \\ \hat{f}_1 \end{Bmatrix} - \delta u_1 \left[ \hat{j}_0 \right] = 0$$

Then, grouping the terms according to the variations and allowing the arbitrary variations on $\delta u_1$, $\delta J_1$, $\delta J_c$, one obtains following equations

$$\left( -\frac{m}{h}\{u_0 - u_1\} + \frac{c}{2}\{u_0 + u_1\} + \frac{1}{6}J_0 + \frac{1}{6}J_1 + \frac{2}{3}J_c - \frac{h}{3}\hat{f}_0 - \frac{h}{6}\hat{f}_1 \right) = \hat{j}_0 \tag{75}$$

$$\frac{a}{3h}\left(7J_0 + J_1 - 8J_c\right) + \frac{5}{6}u_0 + \frac{1}{6}u_1 = 0 \tag{76}$$

$$\frac{a}{3h}\left(-8J_0 - 8J_1 + 16J_c\right) - \frac{2}{3}u_0 + \frac{2}{3}u_1 = 0 \tag{77}$$

Again, with the adoption of the same strategy as Eqs. (46)-(50) in EHP to express $J_c$ in terms of $J_{r-1}$, $J_r$, $u_{r-1}$, and $u_r$, finally we have

$$\begin{bmatrix} \frac{1}{12}\frac{(X + 6cha)}{ha} & \frac{1}{2} \\ \frac{1}{2} & -\frac{a}{h} \end{bmatrix} \begin{Bmatrix} u_1 \\ J_1 \end{Bmatrix} = \begin{bmatrix} \frac{1}{12}\frac{(X - 6cha)}{ha} & -\frac{1}{2} \\ -\frac{1}{2} & -\frac{a}{h} \end{bmatrix} \begin{Bmatrix} u_0 \\ J_0 \end{Bmatrix} + \begin{Bmatrix} Q_{J_1} \\ 0 \end{Bmatrix} \tag{78}$$

where $X$ is defined in Eq. (51) and $Q_{J_1}$ is given by

$$Q_{J_1} = \frac{h}{3}\hat{f}_0 + \frac{h}{6}\hat{f}_1 + \hat{j}_0 \tag{79}$$

More generally, for the $n^{th}$ time step with $t_n = nh$, one may write the Jquad algorithm of MCAP

$$\begin{bmatrix} \frac{1}{12}\frac{(X + 6cha)}{ha} & \frac{1}{2} \\ \frac{1}{2} & -\frac{a}{h} \end{bmatrix} \begin{Bmatrix} u_n \\ J_n \end{Bmatrix} = \begin{bmatrix} \frac{1}{12}\frac{(X - 6cha)}{ha} & -\frac{1}{2} \\ -\frac{1}{2} & -\frac{a}{h} \end{bmatrix} \begin{Bmatrix} u_{n-1} \\ J_{n-1} \end{Bmatrix} + \begin{Bmatrix} Q_{J_n} \\ 0 \end{Bmatrix} \tag{80}$$

where

$$Q_{J_n} = \frac{h}{3}\hat{f}_{n-1} + \frac{h}{6}\hat{f}_n + \hat{j}_{n-1} \tag{81}$$

Similarly, we can develop the Uquad algorithm as

$$\begin{bmatrix} -\frac{1}{3}\frac{(3m+ch)}{h} & \frac{1}{6}\frac{(Y+6cha)}{h^2} \\ \frac{2}{3}\frac{(6m+ch)}{h} & -\frac{2}{3}\frac{(X)}{h^2} \end{bmatrix} \begin{Bmatrix} u_n \\ J_n \end{Bmatrix} = \begin{bmatrix} \frac{1}{3}\frac{(9m+2ch)}{h} & \frac{1}{6}\frac{(Y-6h^2+6cha)}{h^2} \\ \frac{2}{3}\frac{(-6m+ch)}{h} & -\frac{2}{3}\frac{(X)}{h^2} \end{bmatrix} \begin{Bmatrix} u_{n-1} \\ J_{n-1} \end{Bmatrix} + \begin{Bmatrix} Q_{u_n} \\ \frac{h}{3}\hat{f}_{n-1} + \frac{h}{3}\hat{f}_n \end{Bmatrix} \tag{82}$$

where $X$ and $Y$ are given respectively in Eq. (51) and Eq. (53), while $Q_{u_n}$ is given by

$$Q_{u_n} = \frac{h}{6}\hat{f}_{n-1} + \hat{j}_{n-1} \tag{83}$$

Also, we have the UJquad algorithm as

$$\begin{bmatrix} \frac{1}{12}\frac{(mX+6cham+c^2h^2a)}{ham} & \frac{1}{12}\frac{(6m+ch)}{m} \\ \frac{1}{12}\frac{(6m+ch)}{m} & -\frac{1}{12}\frac{(X)}{hm} \end{bmatrix} \begin{Bmatrix} u_n \\ J_n \end{Bmatrix} = \begin{bmatrix} \frac{1}{12}\frac{(mX-6cham+c^2h^2a)}{ham} & \frac{1}{12}\frac{(-6m+ch)}{m} \\ \frac{1}{12}\frac{(-6m+ch)}{m} & -\frac{1}{12}\frac{(X)}{hm} \end{bmatrix} \begin{Bmatrix} u_{n-1} \\ J_{n-1} \end{Bmatrix} + \begin{Bmatrix} Q_{uJ_n} \\ \frac{h^2}{24m}\hat{f}_{n-1} + \frac{h^2}{24m}\hat{f}_n \end{Bmatrix} \tag{84}$$

where

$$Q_{uJ_n} = \left(\frac{h}{3} + \frac{ch^2}{24m}\right)\hat{f}_{n-1} + \left(\frac{h}{6} + \frac{ch^2}{24m}\right)\hat{f}_n + \hat{j}_{n-1} \tag{85}$$

## 4. Basic numerical properties

For the SDOF Kelvin-Voigt model, every algorithm from EHP and MCAP can be written in matrix form as

$$\mathbf{A_1}\,\mathbf{x_n} = \mathbf{A_0}\,\mathbf{x_{n\text{-}1}} + \mathbf{f_n} \tag{86}$$

or simply

$$\mathbf{x_n = A_D\, x_{n-1} + A_1^{-1} f_n} \tag{87}$$

where

$$\mathbf{A_D = A_1^{-1} A_0} \tag{88}$$

### *4.1. Symplectic nature*

For the undamped case with no external forcing (conservative harmonic oscillator), Eqs. (86)-(88) reduce to

$$\mathbf{A_{left}\, x_n = A_{right}\, x_{n-1}} \tag{89}$$

$$\mathbf{x_n = A\, x_{n-1}} \tag{90}$$

$$\mathbf{A = A_{left}^{-1} A_{right}} \tag{91}$$

where $\mathbf{A_{left}}$ and $\mathbf{A_{right}}$ in each algorithm are identified in Table 2 and Table 3.

**Table 2** Algorithms from EHP for the conservative system

| Algorithms | $\mathbf{A_{left}}$ | $\mathbf{A_{right}}$ |
|---|---|---|
| Jquad | $\begin{bmatrix} \frac{X}{12ha} & 0 & \frac{1}{2} \\ -\frac{X}{12ha} & 1 & \frac{1}{2} \\ 1 & 0 & -\frac{2a}{h} \end{bmatrix}$ | $\begin{bmatrix} \frac{X}{12ha} & 1 & \frac{1}{2} \\ -\frac{X}{12ha} & 0 & \frac{1}{2} \\ -1 & 0 & -\frac{2a}{h} \end{bmatrix}$ |
| Uquad | $\begin{bmatrix} -\frac{m}{h} & 0 & \frac{Y}{6h^2} \\ -\frac{3m}{h} & -m & \frac{Y}{6h^2} \\ \frac{4m}{h} & 0 & -\frac{2X}{3h^2} \end{bmatrix}$ | $\begin{bmatrix} \frac{3m}{h} & -m & \frac{Y}{6h^2} \\ \frac{m}{h} & 0 & \frac{Y}{6h^2} \\ -\frac{4m}{h} & 0 & -\frac{2}{3}\frac{X}{h^2} \end{bmatrix}$ |
| UJquad | $\begin{bmatrix} -\frac{1}{12}\frac{\left(12ma+h^2\right)}{ha} & 0 & \frac{Y}{6h^2} \\ \frac{1}{12}\frac{\left(-36ma+h^2\right)}{ha} & 1 & \frac{Y}{6h^2} \\ \frac{4m}{h} & 0 & -\frac{2X}{3h^2} \end{bmatrix}$ | $\begin{bmatrix} \frac{1}{12}\frac{\left(36ma-h^2\right)}{ha} & 1 & \frac{Y}{6h^2} \\ \frac{1}{12}\frac{\left(12ma+h^2\right)}{ha} & 0 & \frac{Y}{6h^2} \\ -\frac{4m}{h} & 0 & -\frac{2X}{3h^2} \end{bmatrix}$ |

**Table 3** Algorithms from MCAP for the conservative system

| Algorithms | $\mathbf{A}_{\mathbf{left}}$ | $\mathbf{A}_{\mathbf{right}}$ |
|---|---|---|
| Jquad | $\begin{bmatrix} \frac{X}{12ha} & \frac{1}{2} \\ \frac{1}{2} & -\frac{a}{h} \end{bmatrix}$ | $\begin{bmatrix} \frac{X}{12ha} & -\frac{1}{2} \\ -\frac{1}{2} & -\frac{a}{h} \end{bmatrix}$ |
| Uquad | $\begin{bmatrix} -\frac{m}{h} & \frac{Y}{6h^2} \\ \frac{4m}{h} & -\frac{2X}{3h^2} \end{bmatrix}$ | $\begin{bmatrix} \frac{3m}{h} & \frac{(Y-6h^2)}{6h^2} \\ -\frac{4m}{h} & -\frac{2X}{3h^2} \end{bmatrix}$ |
| UJquad | $\begin{bmatrix} \frac{X}{12ha} & \frac{1}{2} \\ \frac{1}{2} & -\frac{X}{12hm} \end{bmatrix}$ | $\begin{bmatrix} \frac{X}{12ha} & -\frac{1}{2} \\ -\frac{1}{2} & -\frac{X}{12hm} \end{bmatrix}$ |

In Table 2 and Table 3, $X$ and $Y$ are given respectively in Eq. (51) and Eq. (53), while $Z$ is given by

$$Z = 6ma + h^2 \tag{92}$$

Notice that every algorithm shown in Table 2 and Table 3 is time reversible. One can exactly recover the state $n-1$ from the state $n$ by setting $h \to -h$, $n \to n-1$, and $n-1 \to n$.
For the representative one, one can obtain Uquad algorithm in MCAP as

$$\begin{bmatrix} -\frac{3m}{h} & \frac{Y-6h^2}{6h^2} \\ \frac{4m}{h} & -\frac{2X}{3h^2} \end{bmatrix} \begin{Bmatrix} u_n \\ J_n \end{Bmatrix} = \begin{bmatrix} \frac{m}{h} & \frac{Y}{6h^2} \\ -\frac{4m}{h} & -\frac{2X}{3h^2} \end{bmatrix} \begin{Bmatrix} u_{n-1} \\ J_{n-1} \end{Bmatrix} \tag{93}$$

with the substitution of $h \to -h$, $n \to n-1$, and $n-1 \to n$.
Pre-multiplying the matrix

$$\begin{bmatrix} -1 & -1 \\ 0 & 1 \end{bmatrix} \tag{94}$$

on Eq. (93) yields

$$\begin{bmatrix} -\frac{m}{h} & \frac{Y}{6h^2} \\ \frac{4m}{h} & -\frac{2X}{3h^2} \end{bmatrix} \begin{Bmatrix} u_n \\ J_n \end{Bmatrix} = \begin{bmatrix} \frac{3m}{h} & \frac{(Y-6h^2)}{6h^2} \\ -\frac{4m}{h} & -\frac{2X}{3h^2} \end{bmatrix} \begin{Bmatrix} u_{n-1} \\ J_{n-1} \end{Bmatrix} \tag{95}$$

which is the exactly same Uquad algorithm given in Table 3.

While deriving Eq. (95), the following relation is used

$$-Y + 4X = Y - 6h^2 \tag{96}$$

The stability and dissipative character of each developed method can be determined by considering the eigenvalues of **A** in Eq. (91), and the eigenvalues of each method are presented in Table 4 and Table 5, respectively.

**Table 4** Eigenvalues of A in EHP algorithms

| Algorithms | Eigenvalues |
|---|---|
| Jquad | $\lambda_1 = 1$<br>$\lambda_{2,3} = \dfrac{6ma - 2h^2 \pm i\sqrt{36h^2 ma - 3h^4}}{6ma + h^2}$ |
| Uquad | $\lambda_1 = 1$<br>$\lambda_{2,3} = \dfrac{6ma - 2h^2 \pm i\sqrt{36h^2 ma - 3h^4}}{6ma + h^2}$ |
| UJquad | $\lambda_1 = 1$<br>$\lambda_{2,3} = \dfrac{h^4 - 60amh^2 + 144m^2a^2 \pm i(12h)\left\|12am - h^2\right\|\sqrt{am}}{h^4 + 12amh^2 + 144m^2a^2}$ |

**Table 5** Eigenvalues of A in MCAP algorithms

| Algorithms | Eigenvalues |
|---|---|
| Jquad | $\lambda_{1,2} = \dfrac{6ma - 2h^2 \pm i\sqrt{36h^2 ma - 3h^4}}{6ma + h^2}$ |
| Uquad | $\lambda_{1,2} = \dfrac{6ma - 2h^2 \pm i\sqrt{36h^2 ma - 3h^4}}{6ma + h^2}$ |
| UJquad | $\lambda_{1,2} = \dfrac{h^4 - 60amh^2 + 144m^2a^2 \pm i(12h)\left\|12am - h^2\right\|\sqrt{am}}{h^4 + 12amh^2 + 144m^2a^2}$ |

Notice that the magnitude of all the eigenvalues including complex conjugate pairs in Table 4 and Table 5 is exactly equal to 1, which can be written simply as

$$|\lambda| = 1 \tag{97}$$

Consequently, in addition to being time reversible, all the presented quadratic temporal finite element algorithms are also symplectic, energy conserving, and unconditionally stable for the undamped case.

### *4.2. Period elongation property in each method*

To check the numerical dispersion property in each developed method, the method by [35, 36] is used for free vibration of the undamped oscillator, where the ratio of the time-step $h$ to the natural period $T_n$ is a control parameter. Also, Newmark's constant average acceleration method and Newmark's linear acceleration method are adopted for the references.

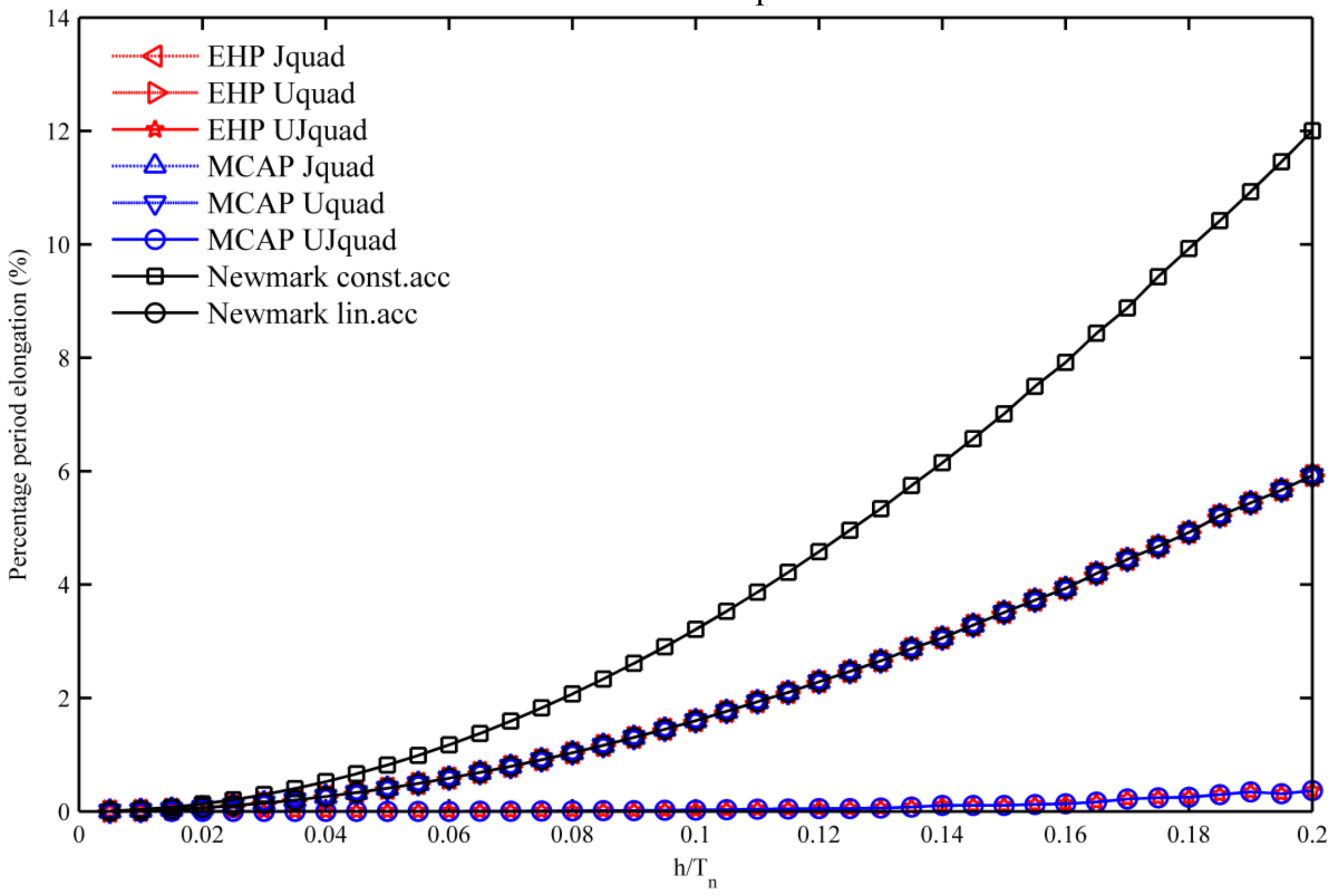

**Fig. 3.** Period elongation property of each method

As shown in Fig. 3, the numerical dispersion property from EHP and MCAP is exactly the same as Newmark's linear acceleration method, when either the primary variable $u$ or $J$ is quadratically approximated. On the other hand, when $u$ and $J$ are quadratically approximated, UJquad algorithm in each method has the same numerical dispersion property better than Newmark's linear acceleration method. Note that all the developed methods are unconditionally stable, while Newmark's linear acceleration method is a conditionally stable algorithm with the criterion

$$\frac{h}{T_n} < 0.551 \tag{98}$$

In computational aspects, compared to Newmark's constant average acceleration and Newmark's linear acceleration method, all the developed computational methods seem practically sufficient and accurate, since they have symplectic, unconditionally stable, and less or equivalent period elongation properties, and this is the main reason that only quadratic temporal finite element methods are developed here.

## 5. Numerical examples

For all of the numerical examples considered here, with no loss of generality, the model parameters are taken in non-dimensional form. In particular, let $m = 1$ and $a = 1/\left(4\pi^2\right)$, thus, providing a natural period $T_n = 1$ in the SDOF Kelvin-Voigt damped oscillator.

Two loading cases with zero initial conditions are considered for numerical simulation. The first one is an applied force in the form $\hat{f}(t) = f_0 \sin(\omega_0 t)$ with $f_0 = 100$ and $\omega_0 = 10$, and the other is 1940 El-Centro loading. The additional parameters for each loading case are summarized in Table 6.

**Table 6** Numerical simulation cases

| Sinusoidal loading $\hat{f}(t) = 100\sin(10t)$ | El-Centro loading |
|---|---|
| $\begin{cases} \text{(i) } h = 0.10 \\ \text{(ii) } h = 0.05 \\ \text{(iii) } h = 0.01 \end{cases}$ while damping coefficient $c = 0.2\pi$ is fixed to deliver a non-dimensional damping ratio $\xi = 0.05$. | $\begin{cases} \text{(i) } \xi = 0.05 \\ \text{(ii) } \xi = 0.03 \\ \text{(iii) } \xi = 0.01 \end{cases}$ while the time step is fixed as $h = 0.02$. |

For the references, the results obtained from each developed method are compared to an exact solution for the sinusoidal loading, while the results from Newmark's linear acceleration method in OpenSees [37, 38] are additionally provided. For El-Centro loading, the results from each developed method are compared to those from Newmark's linear acceleration method in OpenSees.

### *5.1. Simulation results under sinusoidal loading*

Fig. 4 displays the numerical solution of displacement versus time, based upon Newmark's linear acceleration method, while Fig. 5-Fig. 10 are obtained from the developed algorithms.

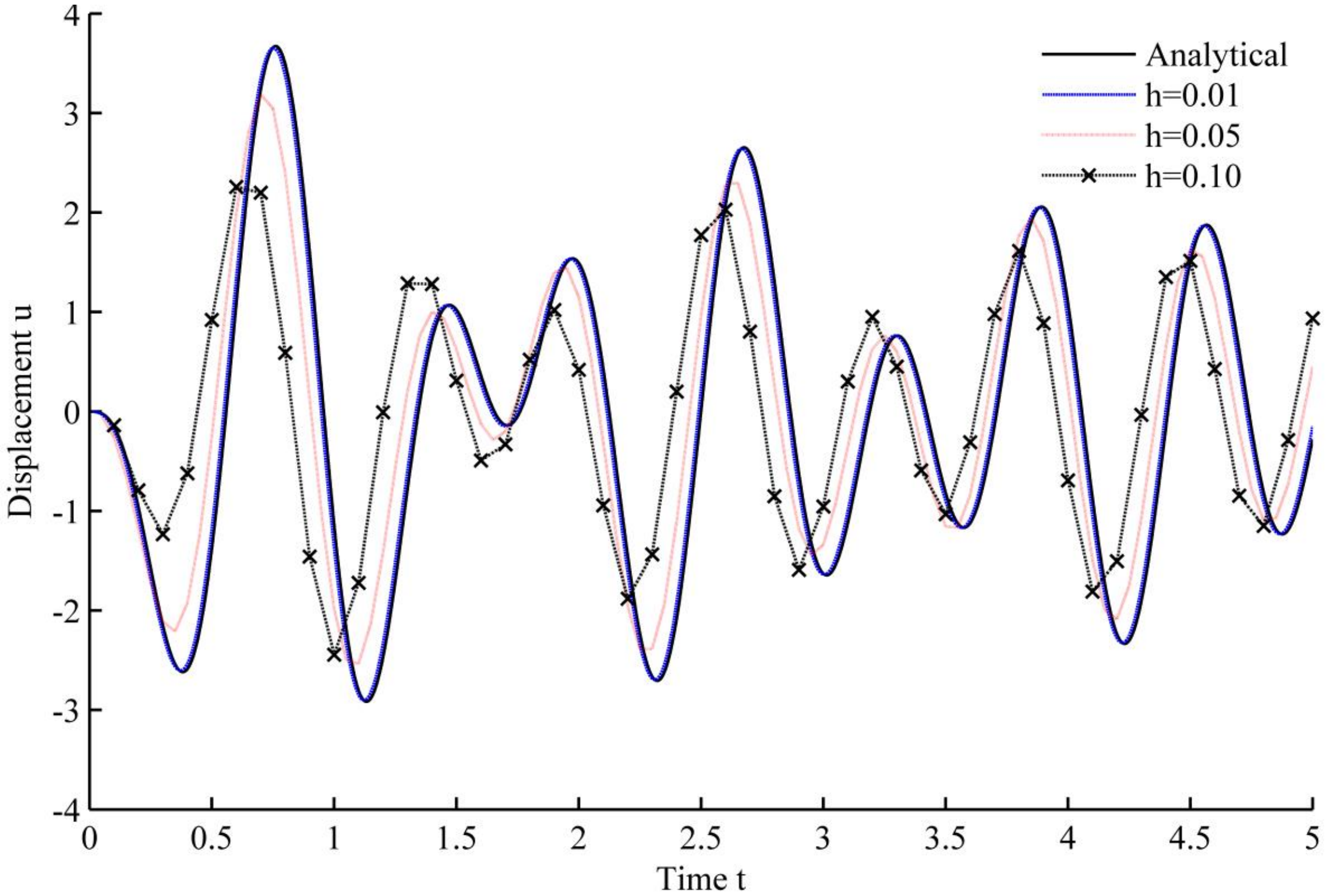

**Fig. 4.** Displacement history results from Newmark's linear acceleration method

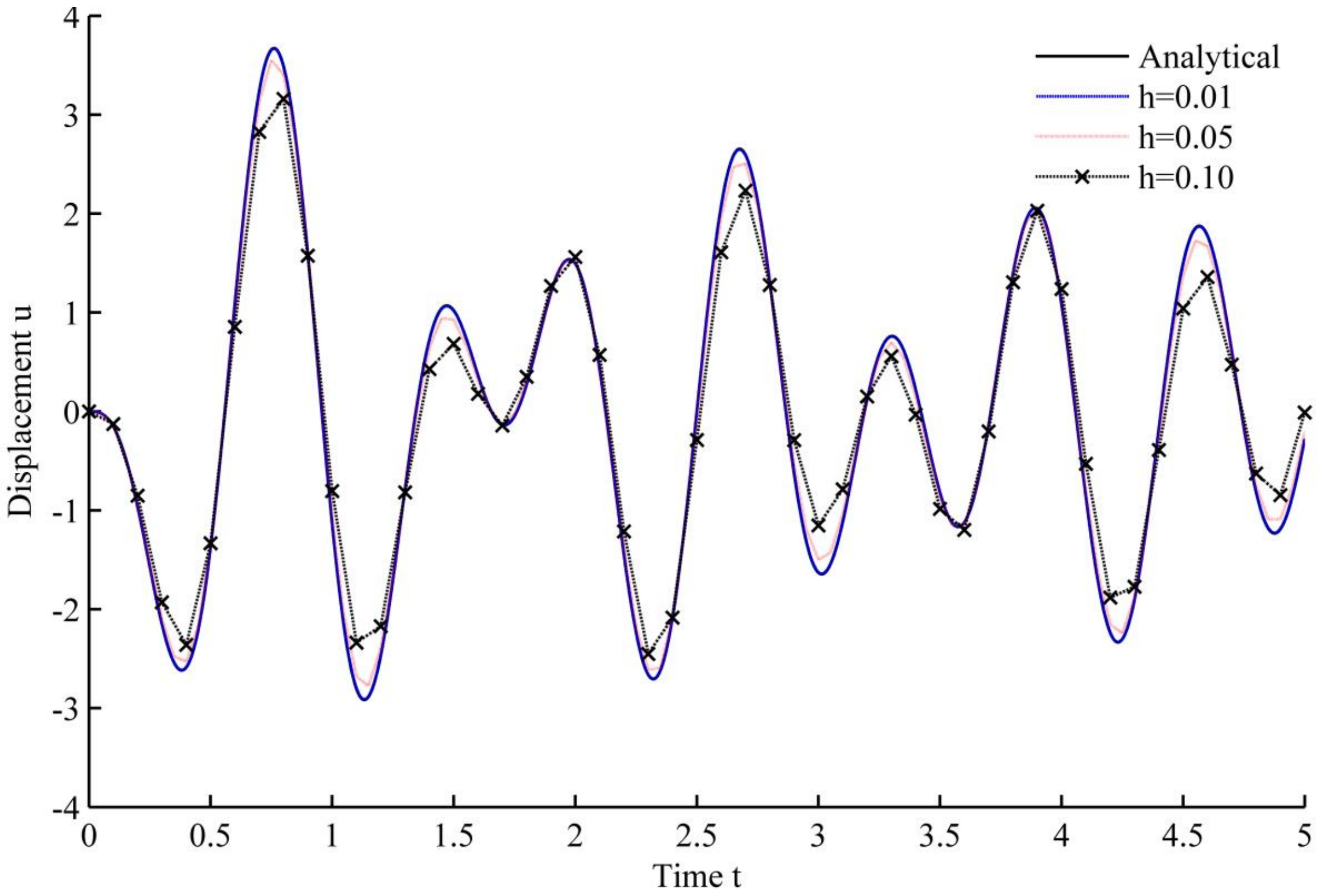

**Fig. 5.** Displacement history results from Jquad algorithm in EHP

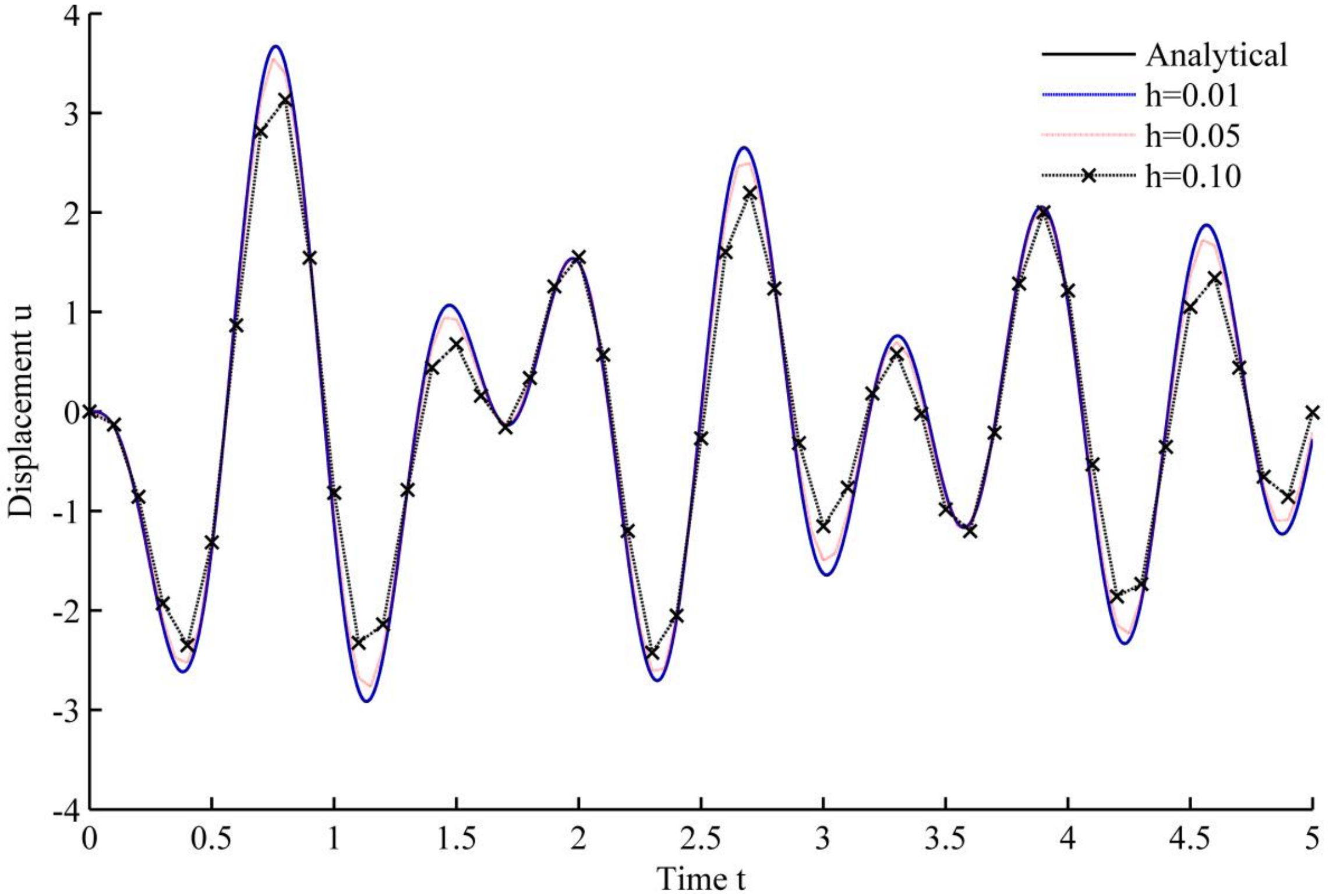

**Fig. 6.** Displacement history results from Uquad algorithm in EHP

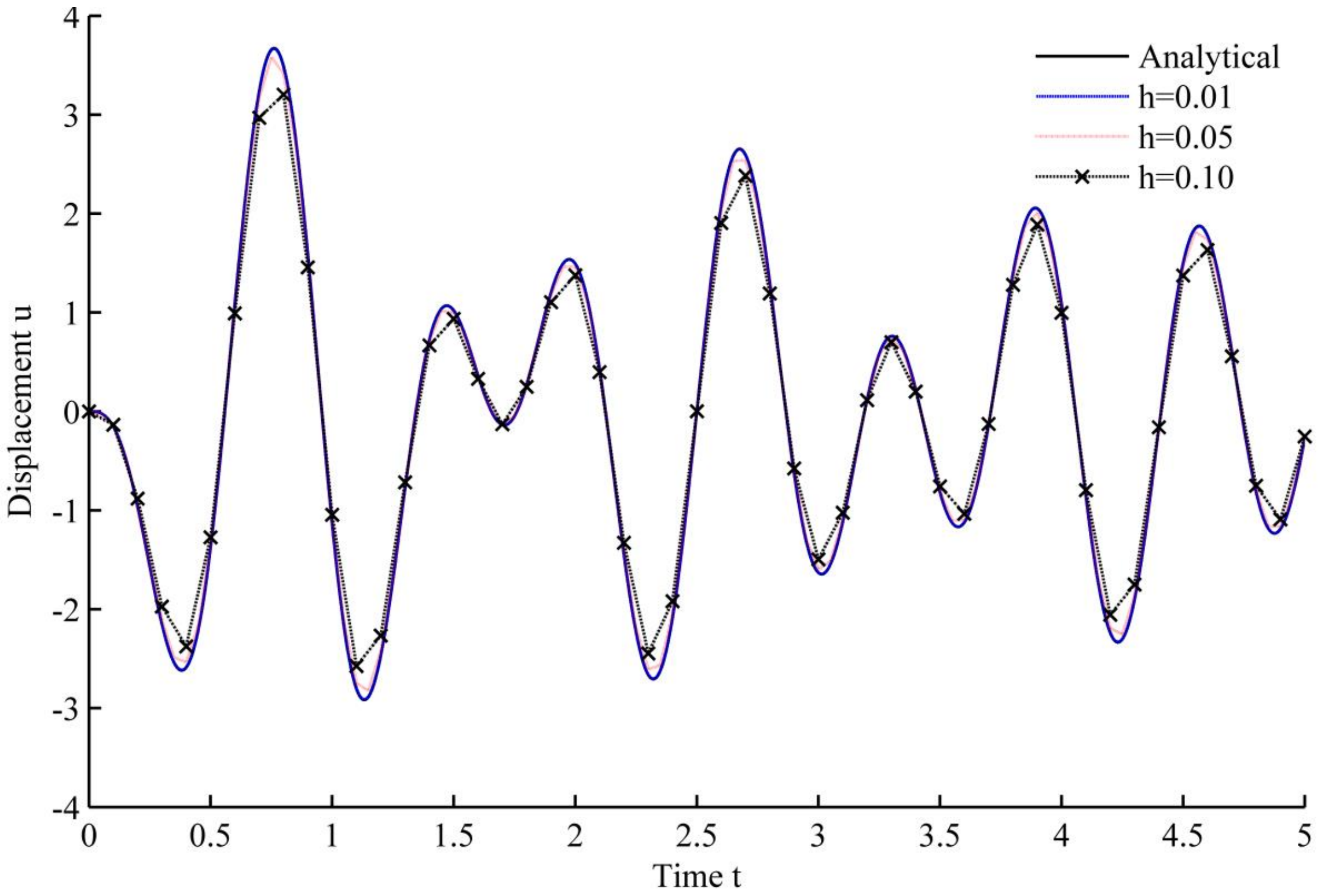

**Fig. 7.** Displacement history results from UJquad algorithm in EHP

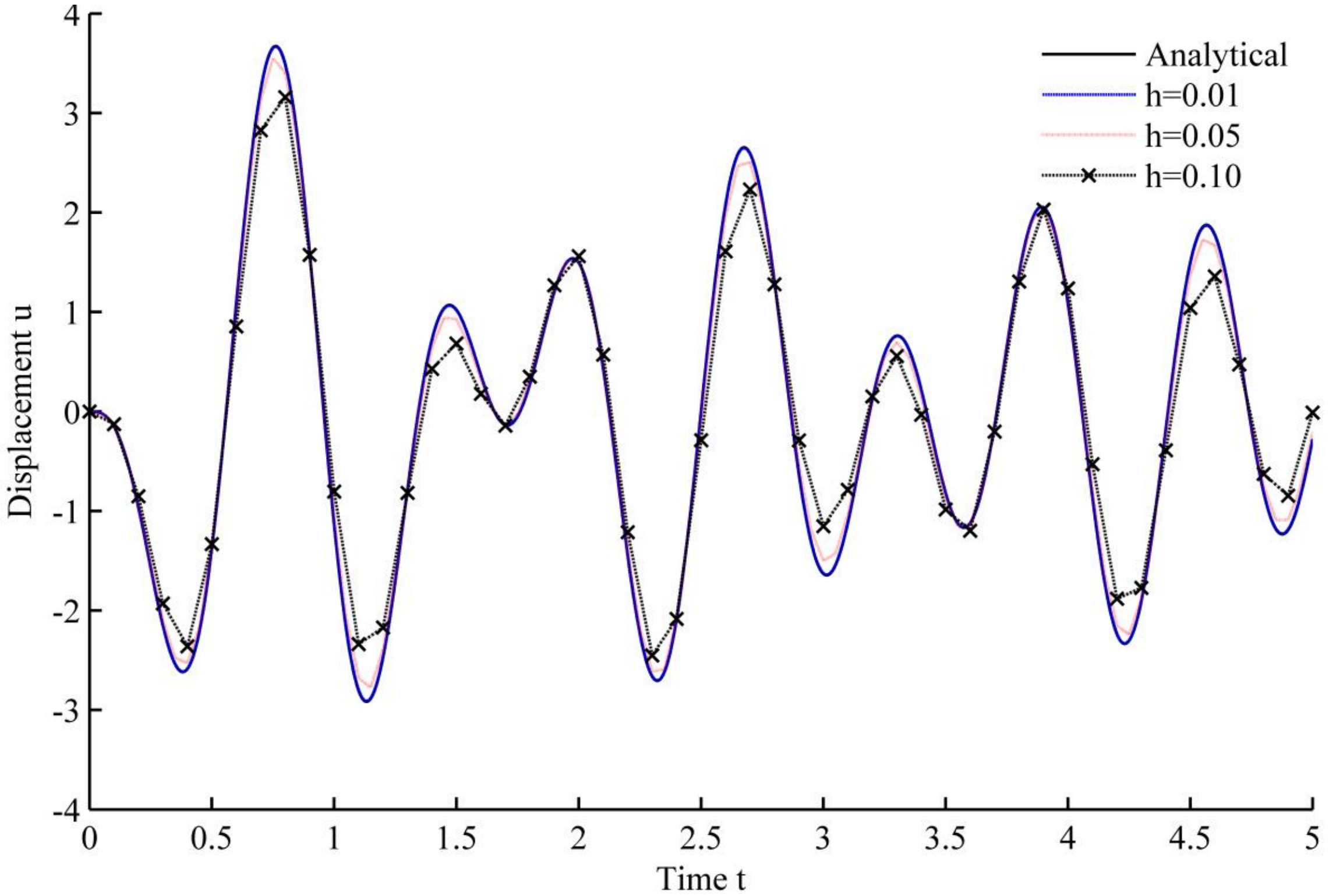

**Fig. 8.** Displacement history results from Jquad algorithm in MCAP

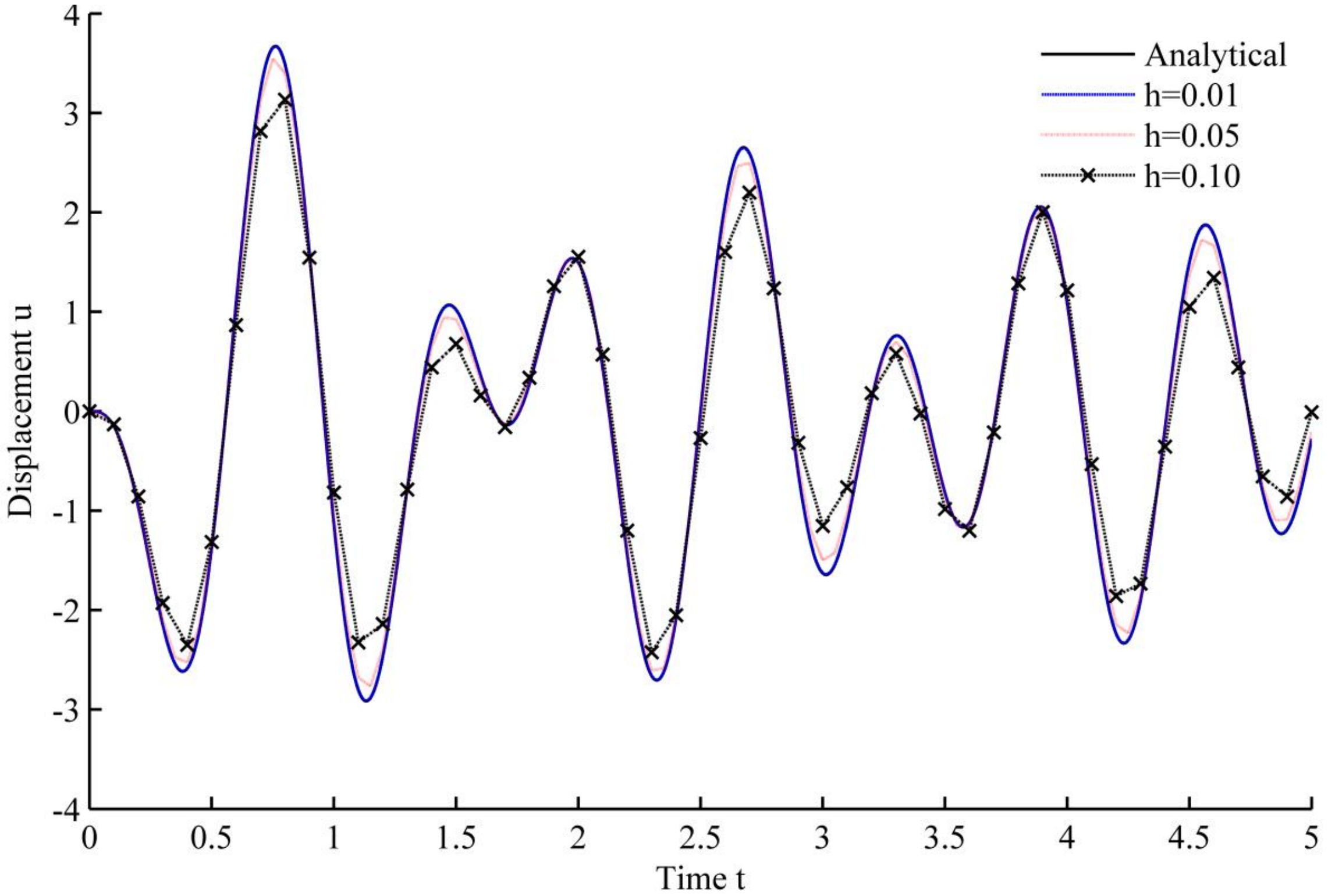

**Fig. 9.** Displacement history results from Uquad algorithm in MCAP

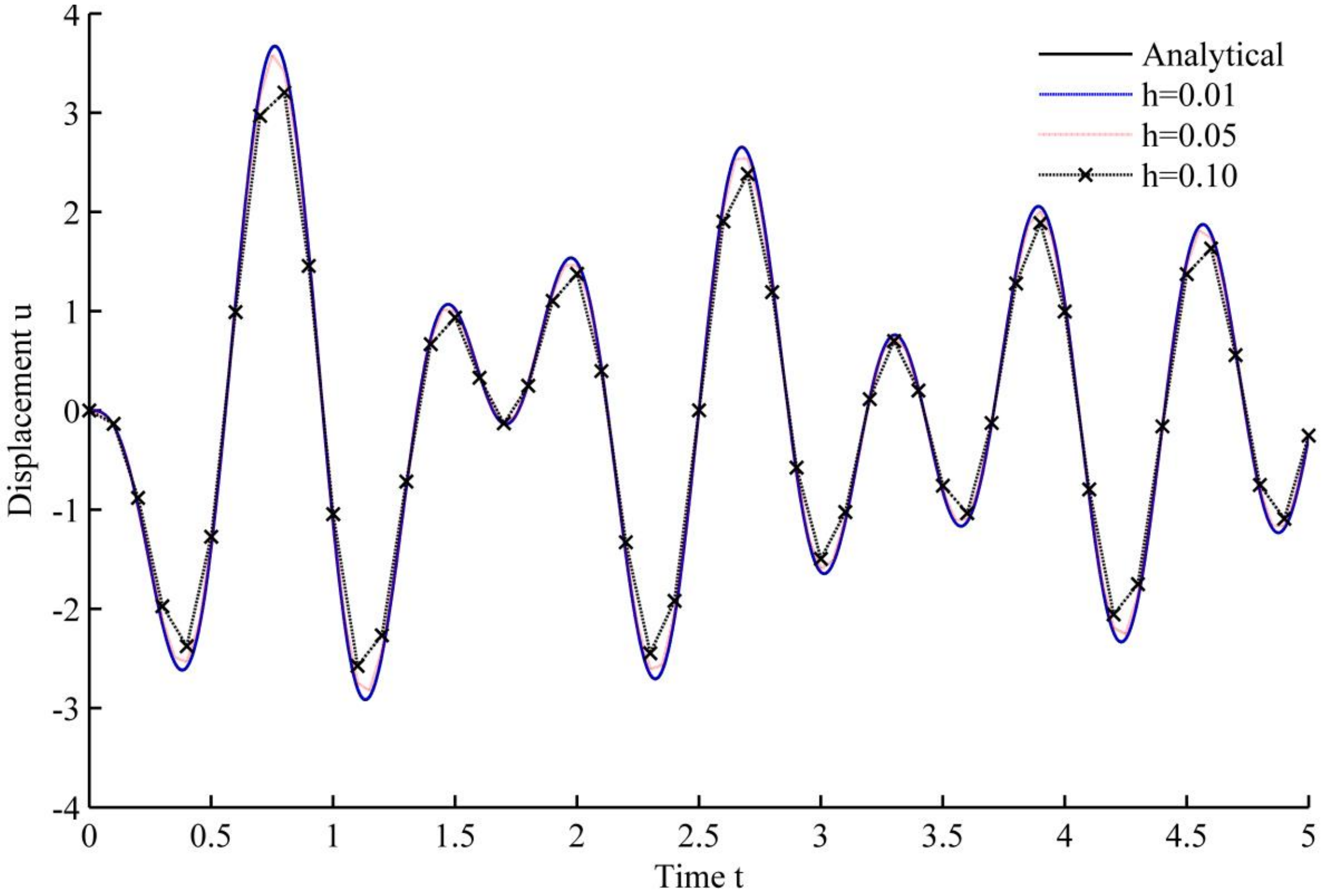

**Fig. 10.** Displacement history results from UJquad algorithm in MCAP

As seen from the results, all the developed methods have better convergence characteristics compared to Newmark's linear acceleration methods under sinusoidal loading. In particular, UJquad algorithm in each framework shows the most accurate results.

### *5.2. Simulation results under 1940 El-Centro loading*

The results from 1940 El-Centro loading analysis are displayed in Fig. 11-Fig. 16. In each figure, the Uquad and Jquad algorithms yield the exactly same results, while there are slight differences between the newly developed methods and Newmark's linear acceleration method. In practical aspects, these differences seem negligible, but, note that all the developed methods are unconditionally stable that it may be advantageous to have the outlined results before the detailed analysis.

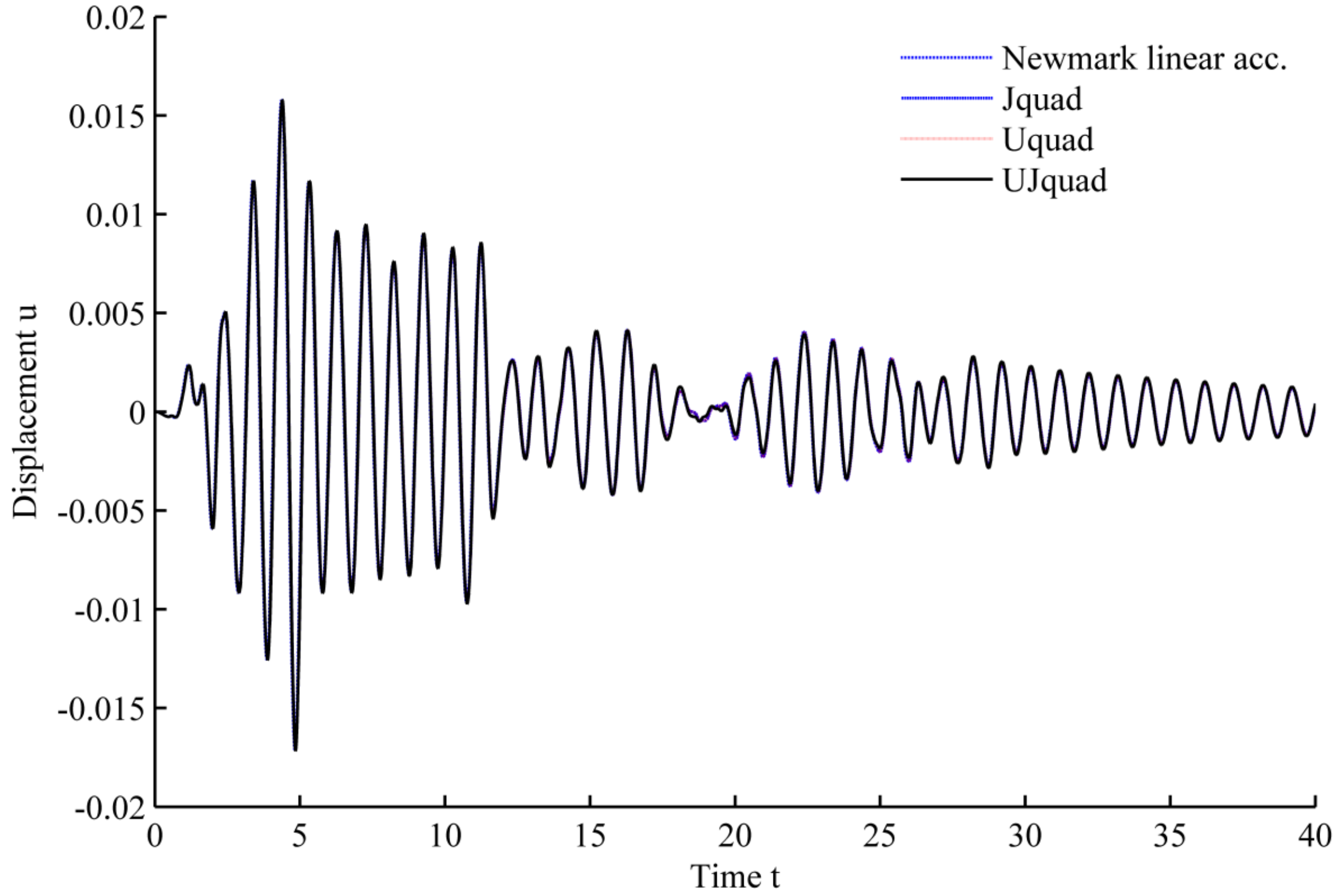

**Fig. 11.** Results from EHP algorithms for El-Centro loading analysis (1% damping ratio)

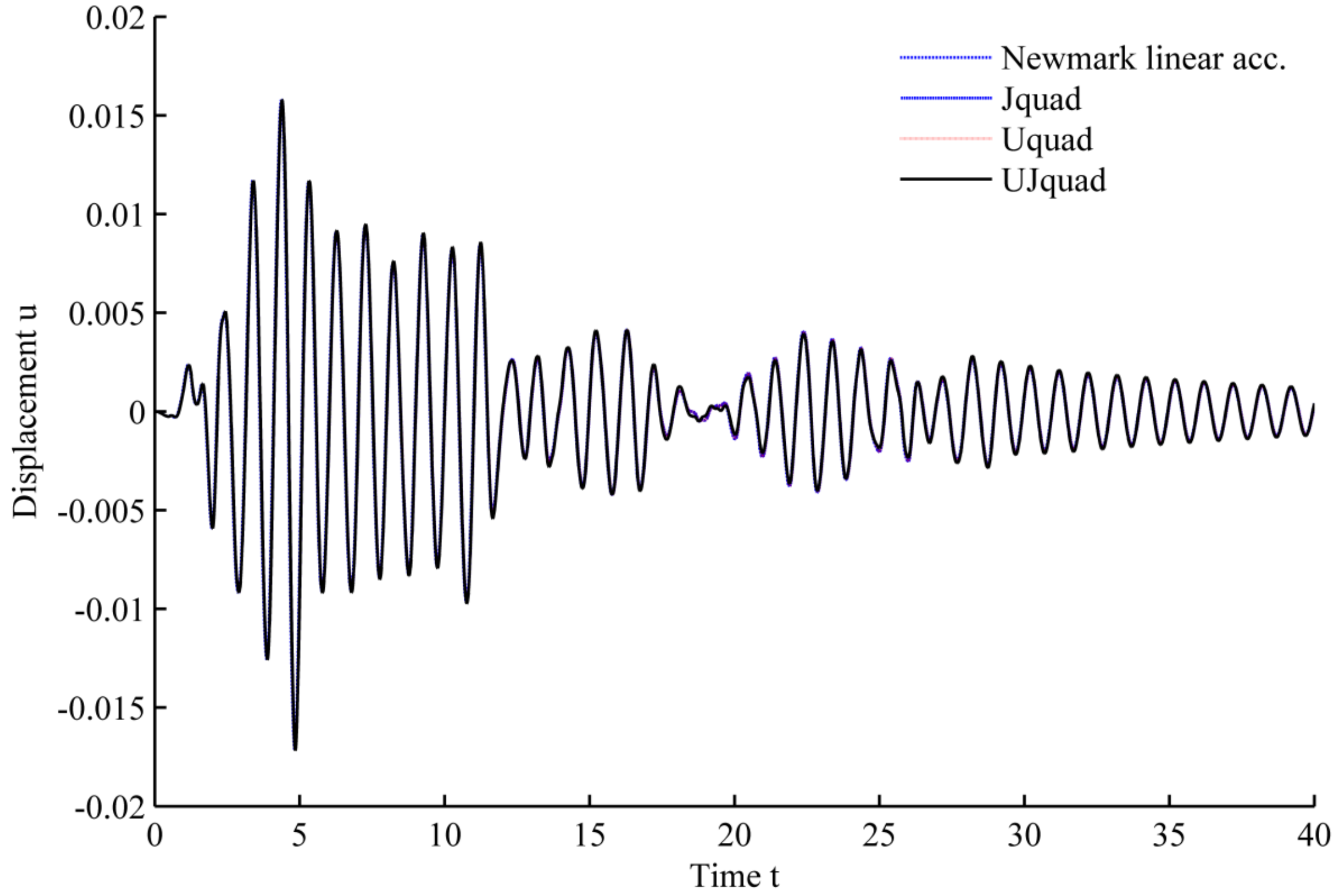

**Fig. 12.** Results from MCAP algorithms for El-Centro loading analysis (1% damping ratio)

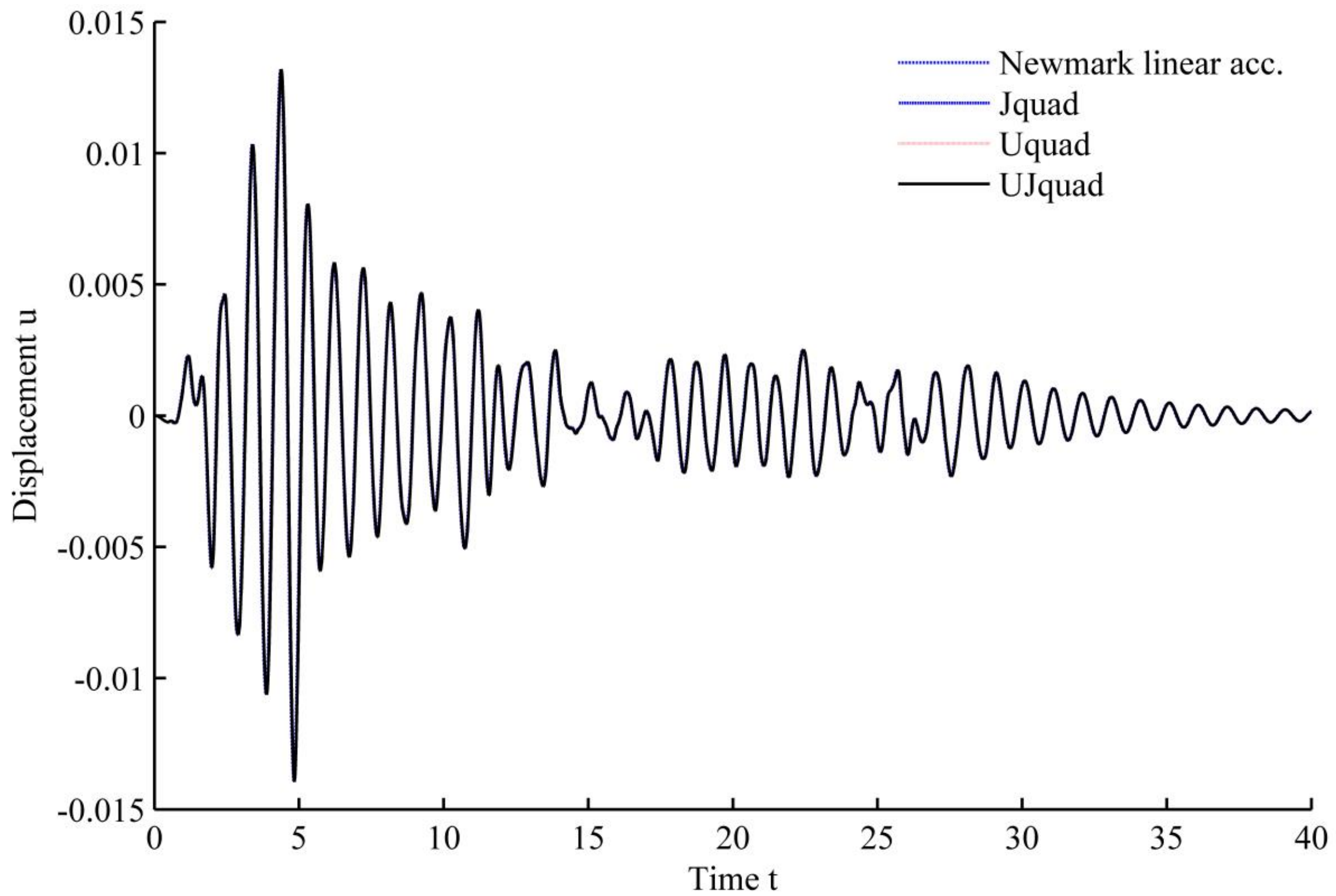

**Fig. 13.** Results from EHP algorithms for El-Centro loading analysis (3% damping ratio)

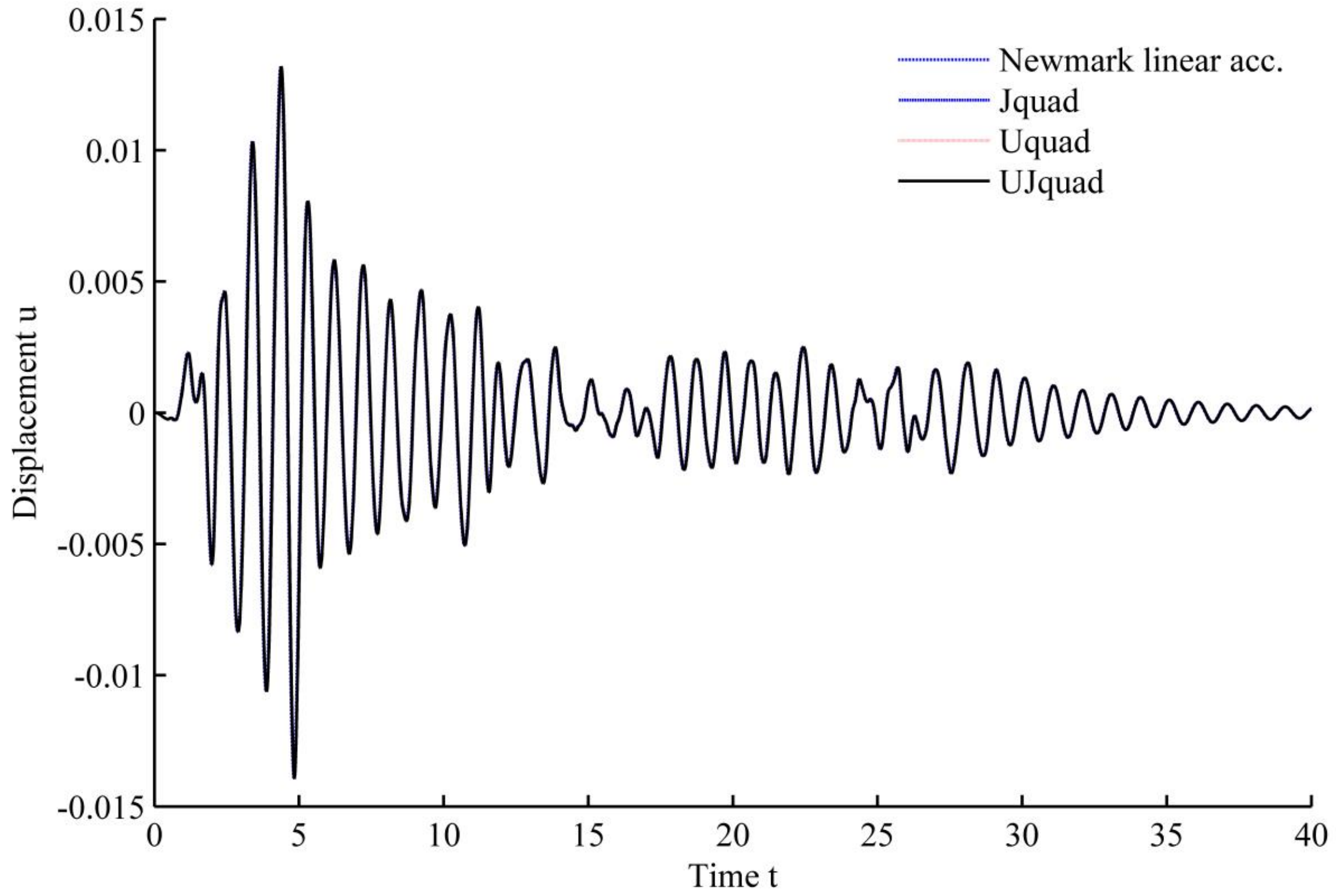

**Fig. 14.** Results from MCAP algorithms for El-Centro loading analysis (3% damping ratio)

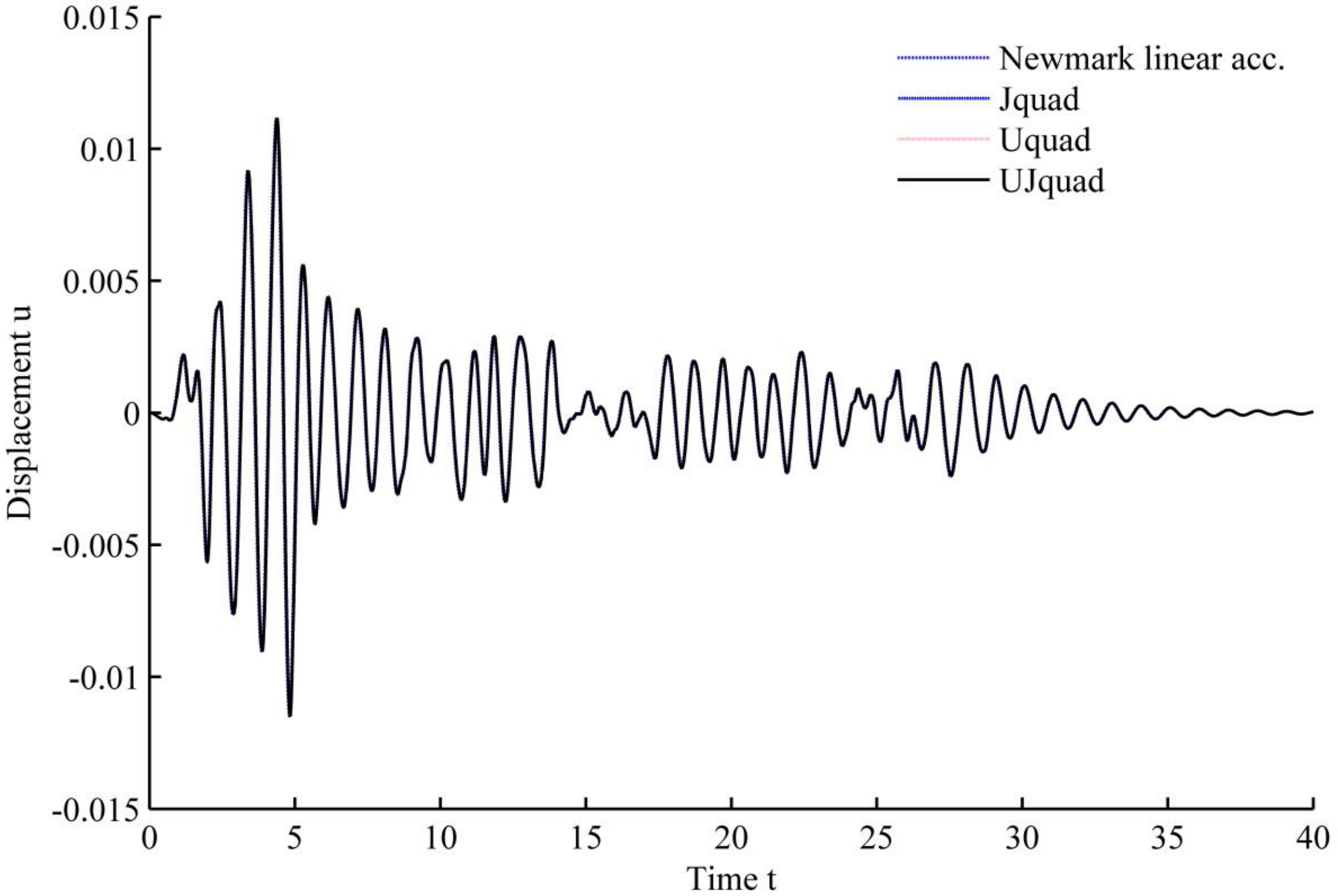

**Fig. 15.** Results from EHP algorithms for El-Centro loading analysis (5% damping ratio)

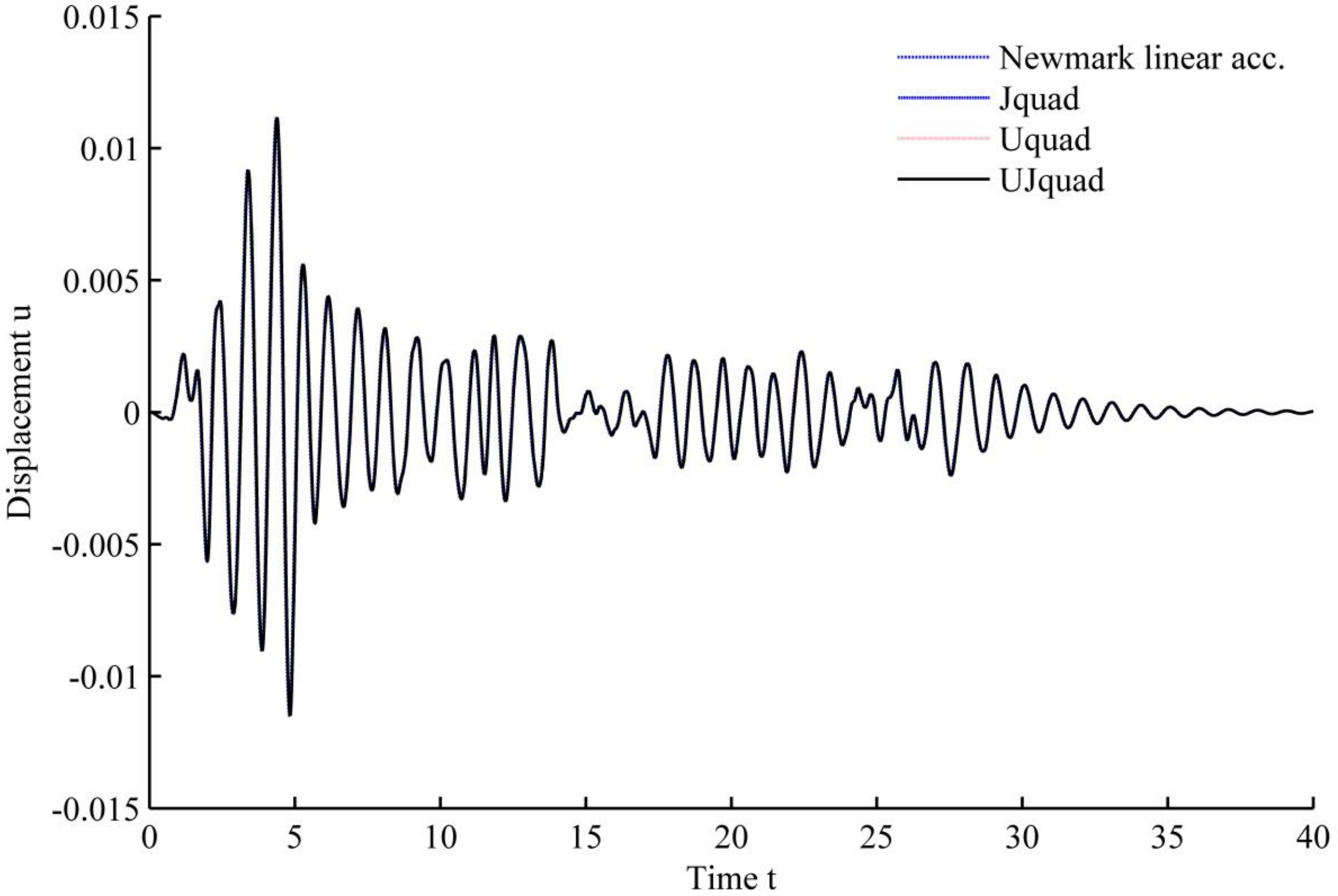

**Fig. 16.** Results from MCAP algorithms for El-Centro loading analysis (5% damping ratio)

## 6. Conclusions

In recent papers, through mixed formulation, two new variational frameworks such as EHP and MCAP were formulated for dynamical systems. Theoretically, MCAP is preferred to EHP, because unlike previous variational approaches, MCAP does not require any dissipation function with ad-hoc rules for taking variations, restrictions on the variations at the ends of the time interval, and external specification of initial conditions. However, there still remains a challenge for MCAP to have a generalized framework embracing various irreversible phenomena. On the other hand, EHP has a relatively simple framework: the action variation is newly defined by adding the counterparts to the terms without the end-point constraints in Hamilton’s principle, which confines a dynamical system to evolve uniquely from start to end. Interpreting these additional terms as sequentially assigning the known initial values completes this formulation. It should be noted that EHP is not a complete variational method, since it still requires the Rayleigh’s dissipation for a non-conservative process and it cannot define the functional action explicitly. Since both mixed formalism provide a rigorous foundation to develop various temporal finite element methods for linear elasticity, in this paper, their potential when adopting temporally higher order approximations is investigated for the classical SDOF Kelvin-Voigt damped system.

With the consideration of computational aspects, three quadratic temporal finite element methods are essentially developed from each mixed formalism. All the developed methods are symplectic and unconditionally stable for the undamped conservative harmonic oscillator. Also, from period elongation property studies, it is checked that all the developed methods are equivalent or superior to Newmark’s linear acceleration method that is conditionally stable. For damped forced vibrations, all the developed methods are shown to be robust and to be accurate with good convergence characteristics. It should be noted that since the new methods utilize mixed formulations, there exists an inherent disadvantage in a significant increase of the degrees of freedom against Newmark’s methods when dealing with other than SDOF systems. However, this may be somewhat compensated by the general characteristics of a mixed formulation and its broad applicability [39-42]

As the original Hamilton’s principle has been adopted in various applications, the applicability of EHP and MCAP are quite broad, spanning many fields of mathematical physics and engineering. Future work will be directed toward development of a generalized framework of MCAP, and applications of both formalisms to various engineering problems, following the ideas in [43-46].